\documentclass[11pt,a4paper]{amsart}
\usepackage{amsfonts}
\usepackage{amsfonts}
\usepackage{amsmath,xypic}
\usepackage{mathrsfs}
\usepackage{amssymb}

\usepackage{CJK,CJKnumb}
\usepackage[CJKbookmarks,colorlinks,
            linkcolor=black,
            anchorcolor=black,
            citecolor=blue]{hyperref}
\usepackage{color}
\usepackage{indentfirst}
\usepackage{latexsym,bm}
\usepackage{amsmath,amssymb}
\usepackage{pstricks}
\usepackage{pst-node}
\usepackage{pst-tree}
\usepackage{pst-plot}
\usepackage{pst-text}
\usepackage{graphicx}
\usepackage{cases}
\usepackage{pifont}
\usepackage{txfonts}
\usepackage[all,knot,poly]{xy}
\allowdisplaybreaks[4]

\newcommand{\R}{\mathscr{R}}

\newtheorem{theorem}{Theorem}[section]
\newtheorem{lemma}{Lemma}[section]
\newtheorem{remark}{Remark}[section]
\newtheorem{proposition}{Proposition}[section]

\allowdisplaybreaks[4]

\begin{document}
\title[Grafting of quantum groups]{Double-bosonization and Majid's conjecture \\(V): grafting of quantum groups}
\author[H. Hu]{Hongmei Hu}
\address{Department of Mathematics, Shanghai Maritime University, Shanghai 201306, PR China}
\email{hmhu@shmtu.edu.cn}

\author[N. Hu]{Naihong Hu$^{\ast}$}
\address{School of Mathematical Sciences, MOE Key Laboratory of Mathematics and Engineering Applications \& Shanghai Key Laboratory of Pure Mathematics and Mathematical Practice,
East China Normal University,
Minhang Campus,
Dong Chuan Road 500,
Shanghai 200241,
PR China}
\email{nhhu@math.ecnu.edu.cn}
\date{\today}
\thanks{$^{\ast}$ The corresponding author.}

\date{}
\maketitle

\newcommand*{\abstractb}[3]{ %
                             \begingroup%
                             \leftskip=8mm \rightskip=8mm
                             \fontsize{11pt}{\baselineskip}\noindent{\textbf{Abstract} ~}{#1}\\
                             {\textbf{Keywords} ~}{#2}\\
                             {\textbf{MR(2010) Subject Classification} ~}{#3}\\
                                                          \endgroup
                           }
\begin{abstract}
This paper aims to develop a grafting method to address Majid’s conjecture. This method enables the construction of a larger target quantum group by grafting two given smaller quantum groups, and
advances the study of the generation, classification, and construction of (quasi‑)Hopf algebras. As a foundation for this construction, we establish a multi-tensor-product theory of generalized double-bosonization and use it to extract the crucial data for the associated braiding $R$-matrix.
Beyond the braided monoidal category perspective arising from quantum subgroup representations, the grafting procedure needs to incorporate structural information from Lie-theoretic root systems. The resulting theoretic framework provides a one-stop strategy for resolving the generation problem concerning the quantum-groups tree in Majid’s conjecture.
\end{abstract}
\section{Introduction}

\noindent{\bf 1.1.} This paper continues our series of papers \cite{HH1,HH2,HH3} on Majid’s conjecture, which
asserts that every quantum group of Drinfeld-Jimbo type can be constructed from $U_q(\mathfrak{sl}_2)$ by a sequence of so-called double-bosonization procedures.
Each such double-bosonization consists of two Radford biproducts \cite {S} --- that is,
bosonizations of braided groups in Majid’s sense \cite {M,majid1,majid2,majid3,Gra} --- that
share a common Cartan part consisting of group-like elements, as in Radford’s projection
construction \cite {rad,rad1}.
This double-bosonization construction, which is closely related to the quantum-double
structure, has been extensively generalized within the theory of quasi-Hopf algebras (or
quasi-quantum groups) initiated by Drinfeld \cite {D}; it continues to play a central role in
current research (cf.\ \cite {BCPO,BN,BPT,BT,BT1,BT2}). Its dual version,
co-double-bosonization \cite {AM}, has likewise influenced recent developments \cite
{AAY,BC,B}.

It is well known that studying the structural theory of (quasi-)quantum groups within the framework of (quasi-)Hopf algebras is more appropriate, as such, this enables the identification of additional distinctive features that set them apart from those in Lie theory (\cite{BD, M, majid4, M5, S}).
  The theory of Hopf algebras has developed a profound body of theories, such as braided tensor categories (\cite{HS}, \cite{BCPO}), which are closely related to topological quantum field theory (\cite{TV}, etc.). Over the past quarter-century, the classification of finite-dimensional pointed Hopf algebras over the complex field by Andruskiewitsch, Heckenberger, Schneider (\cite{AS, HS}) is a major milestone in Hopf algebra theory. This work is closely linked to the classification of braided groups, also known as Nichols algebras (\cite{HS}), in fact, recently it has been found that has a close relationship with the structure theory of multi-parameter quantum groups at roots of unity (\cite{BGM}, \cite{GG}, \cite{PHR}, etc.).
In classification theory, the construction or reconstruction of Hopf algebras (\cite{majid1, majid2, majid3, ro, So}) is crucial in a certain sense.
The reconstruction method developed here for standard quantum groups when $q$ is not a root of
unity may be extended either to small multiparameter quantum groups at roots of unity, in
connection with \cite {CL}, or to finite quasi-quantum groups \cite {B,BT1,BT2, HLYY}.
These extensions are challenging problems closely intertwined with the corresponding
representation theories and may lead to further structural developments.

\medskip
\noindent{\bf 1.2.}
Since Majid reformulated the Radford-Majid bosonization in the setting of a braided category over
a quasitriangular Hopf algebra \cite {majid1,majid2,rad}, its indispensable role in the
classification of Hopf algebras and quantum groups has become increasingly clear.
For instance, starting from dual braided groups \(B^\star\) and $B$ in the braided tensor category consisting of $H$-modules for a quasitriangular Hopf algebra $H$, one gets a new quantum group \(B^\star\otimes H\otimes B\); see Majid’s double-bosonization theorem in \cite{majid3}, where the bosonizations \(H\ltimes B\) and \(B^\star\rtimes H\) appear as the positive and negative Borel subalgebras, respectively.
Majid conjectured in \cite{majid2} that \(B^\star\otimes\tilde H\otimes B\) is a quantum enveloping algebra of rank one higher when $\widetilde H$ is a suitable central
extension of $U_q(\mathfrak g)$ for the $ABCD$ series.
Moreover, he suggested in \cite{majid3} that more general considerations for inductive constructions are possible, and that a tree of quantum groups arising from double-bosonization theory may grow starting from \(U_q(\mathfrak{sl}_2)\).
This program was carried out in \cite {HH1,HH2,HH3} for types $A$, $B$, $C$, $D$, and $G_2$.
Figure $1$ illustrates the resulting growth tree for finite-type quantum groups, rooted at \(U_q(\mathfrak{sl}_2)\).
Originally, the tree grows by successively appending a new simple root to one end of a Dynkin diagram. We now ask whether this quantum tree admits
alternative graft-type growth modes: specifically, how one may graft two or finitely many smaller quantum groups to construct a larger one.
We call this operation the grafting of quantum groups. This is an effective approach for yielding new quantum groups; additional examples were announced in \cite{HHX}.

\medskip
\noindent{\bf 1.3.}
The inverse procedure of grafting consists in deleting an interior vertex from the Dynkin diagram of a Drinfeld-Jimbo quantum group.
To formulate the grafting construction, we first establish several standard facts about the data involved in the tensor product of finitely many quantum
groups.  These are necessary prerequisites for setting up the grafting method.
We then identify the key ingredients for grafting two quantum groups. Given two quantum groups, we connect their respective Dynkin diagrams by inserting a middle black dot, thereby obtaining the Dynkin diagram of the resulting larger quantum group. Here, if the connecting line is simply-laced, our multi-tensor variant of double-bosonization theory applies directly. For example, in the type $A$ case shown in Figure~2, $U_q(\mathfrak {sl}_{m+n})$ is constructed
directly by grafting $U_q(\mathfrak {sl}_n)$ onto $U_q(\mathfrak {sl}_m)$, rather than by
applying the rank-inductive construction of \cite {HH1} successively $m-1$ times to
$U_q(\mathfrak {sl}_n)$.
When a multiple lace appears attached to the middle black vertex to one of the two sub-Dynkin diagrams, the
construction becomes more involved than in the $G_2$ case \cite {HH2}.
The higher $q$-Serre relations involving the connecting vertex are encoded within the radicals of the pairings for braided (co‑)vector algebras in Majid’s sense \cite{majid1}. These radicals are determined by the tensor product of the $R$-matrices defined on the respective chosen modules of the two given quantum subgroups.
$$\begin{array}{ccc}
\includegraphics[scale=0.4]{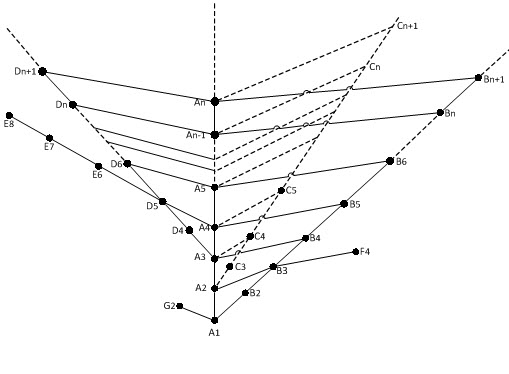}&&
\setlength{\unitlength}{1mm}
\begin{picture}(10,40)
	\put(0,10){$A_{m-1}$}
	\put(0,32){$A_{n-1}$}
	\put(9.5,0){\circle{1}}
	\put(9.5,1){\line(0,1){3}}
	\multiput(9.5,4)(0,1){3}{\line(0,1){0.5}}
	\put(9.5,7){\line(0,1){3}}
	\put(9.5,10.5){\circle{1}}
	\put(9.5,11){\line(0,1){3}}
	\put(9.5,15){\circle{1}}
	\multiput(9.5,15.5)(0,1){3}{\line(0,1){0.5}}
	\put(9.5,19.5){\circle*{1}}
	\multiput(9.5,20.5)(0,1){3}{\line(0,1){0.5}}
	\put(9.5,24){\circle{1}}
	\put(9.5,24.5){\line(0,1){3}}
	\multiput(9.5,28)(0,1){3}{\line(0,1){0.5}}
		\put(9.5,31){\line(0,1){3}}
\put(9.5,34.5){\circle{1}}	
\put(9.5,35.5){\line(0,1){3}}
\put(9.5,39){\circle{1}}	
\end{picture}
\\
\mbox{Figure 1: Tree of quantum groups}\smallskip&&
\qquad\qquad\mbox{Figure 2: Type $A$}
\end{array}$$
This requires taking quotients in the double-bosonization construction by the radicals of the
braided dual pairing. As a representative example, we treat the more complicated case of type
$F_4$; see Figure~3 in Subsection 4.2.

\medskip
\noindent{\bf 1.4.}
Recall that to produce the exceptional $G_2$ branch of the quantum tree, we introduced the generalized double-bosonization construction in \cite{HH2}. This extends Majid’s double-bosonization construction, where the latter applies only to $ABCD$-branches of the quantum tree (\cite{HH1}).  For our purposes here, we need to lay out the foundational background for multi-tensor-products within generalized double-bosonization. This enables us to graft two or finitely many quantum groups equipped with their respective modules and thereby construct the desired larger quantum group.
 As illustrative examples, in the simply laced setting, we explicitly present the grafting construction yielding \(U_q(\mathfrak{sl}_{n+m})\) from the pair \(\bigl(U_q(\mathfrak{sl}_n),\,U_q(\mathfrak{sl}_m)\bigr)\) using their (dual) natural modules,
and in the non-simply laced setting, we construct $U_q(F_4)$ via grafting $U_q(\mathfrak {sl}_3)$ equipped with its natural module onto $U_q(\mathfrak{sl}_2)$ with its natural module.
We emphasize here that the higher q-Serre relations, as involved in \(U_q(F_4)\), arise by quotienting the enlarged braided (co)vector algebras modulo the radical of their pairing.
Indeed, such a grafting method yields further examples and even new quantum groups beyond finite types, including affine and even indefinite types. An example concerning the strictly hyperbolic Kac-Moody algebra of type \(T_{p,q,r}\) (cf. Exercise 4.2 in \cite{K}) was presented in \cite{HHX}.
Formulating a quantum analogue of the pivotal work by Benkart-Kang-Misra (\cite{BKM}) for infinite-dimensional classical-type Kac-Moody algebras constitutes a challenging task for future work.

\noindent{\bf 1.5.}
The paper is organized as follows. We collect some basic facts on generalized double-bosonization theory in Section 2.
In Section~3, we first derive an explicit formula for every entry of the $R$-matrix
$R_{V_1\otimes \cdots \otimes V_n}$ from the universal $R$-matrices of the factors; see
Proposition~\ref {rmatrix}. Proposition~\ref {minimal} then determines the eigenvalues, and
hence the minimal polynomial, of the associated braiding matrix $PR_{V_1\otimes \cdots \otimes
V_n}$, thereby yielding a Majid pair $(R,R’)$ satisfying \eqref {*} and \eqref {**}.
Second, recall that Majid’s double-bosonization theory requires one to find dually paired braided groups \(B^\star\) and $B$ lying in the braided category of $H$-modules, such that the tensor space \(B^\star\otimes H\otimes B\) carries the structure of a new Hopf algebra. To achieve this, we require a weakly quasitriangular dual pair \((H,A)\), where $A$ is a coquasitriangular Hopf algebra.
Choosing $U_q(\mathfrak{g}_i)$ as $H_i$ for some semisimple Lie algebra $\mathfrak g_i$,
Propositions~\ref {hopf2} and \ref {pm} show that the coquasitriangular Hopf subalgebra
$H_{R_{V_1\otimes \cdots \otimes V_n}}$ is the appropriate choice for $A$.
We prove in Theorem \ref{weakly22} that  $U_q^{\textrm{ext}}(\mathfrak{g}_1\oplus\cdots\oplus\mathfrak{g}_n)$ and $H_{R_{V_1\otimes\cdots\otimes V_n}}$ indeed constitute a weakly quasitriangular dual pair. These results yield
the tensor-product version of generalized double-bosonization stated in Theorem~\ref {multi}.
As applications of Theorem \ref{multi},
when performing the grafting of two quantum groups in Section 4,
the main tasks are to derive the required algebraic relations and to construct bases of braided groups  $B^{\star}, B$
such that the Hopf algebra constructed on the tensor space $B^{\star}\otimes (\widetilde{U_{q}^{\textrm{ext}}}(\mathfrak{g}_1\oplus\cdots\oplus\mathfrak{g}_n))\otimes B$ is isomorphic to the desired quantum enveloping algebra.
To this end,
we appeal to the Diamond Lemma in \cite{Ber}.
As an example  of simply-laced case,
we construct the quantum group $U_q(\mathfrak{sl}_{m+n})$ by grafting $U_q(\mathfrak{sl}_n)$ and $U_q(\mathfrak{sl}_m)$ with the respective modules $\mathbb{C}^n$ and $(\mathbb{C}^m)^*$ in Theorem \ref{typeA}.
Here,
$\mathbb{C}^n$ is the natural module of $U_q(\mathfrak{sl}_n)$,
and $(\mathbb{C}^m)^*$ is the dual natural module of $U_q(\mathfrak{sl}_m)$.
As an example of nonsimply laced case,  we construct $U_q(F_4)$ on the tensor-product vector space
 $\bar V^{\vee}(R^{\prime},R_{21}^{-1})\otimes \widetilde {U_q^{\textrm{ext}}}(\mathfrak{sl}_3\oplus \mathfrak{sl}_2)\otimes \bar V(R^{\prime},R)$ in Theorem \ref{typeF}.
Here,
the Majid pair $(R, R')$ comes from the tensor product of $R$-matrices of the natural representations of $U_q(\mathfrak{sl}_3)$ and $U_q(\mathfrak{sl}_2)$.
The defining relations, bases of the two braided groups, and the remaining details are given
in Section~4.

\section{Generalities on Double-Bosonization Theory}
\subsection{Some basic facts and notation}\label{fact}
Here, we collect basic facts and notation on braided groups, the FRT construction, and double-bosonization that are needed in the following sections.
For a more detailed review,
see \cite{FRT,majid1,HH2}.

\subsubsection{Majid's weakly quasitriangular dual pair}\label{gel}
Majid \cite{majid2} proposed the concept of a weakly quasitriangular dual pair via his insight on some examples on
matched pairs of bialgebras or Hopf algebras.
\emph{	A weakly quasitriangular dual pair} is a pair $(H, A)$ of bialgebras or Hopf algebras equipped with a dual pairing
	$\langle \,,\,\rangle: H\otimes A \rightarrow k$
	and convolution-invertible algebra\,/\,anti-coalgebra maps $\mathcal{R},\bar{\mathcal{R}}: A \rightarrow H$ obeying
	\begin{equation}\label{weak:pair}
		\langle\bar{\mathcal{R}}(a),b\rangle=\langle\mathcal{R}^{-1}(b),a\rangle,\quad
		\partial^{R}_h=\mathcal{R} \ast(\partial^{L}_h)\ast\mathcal{R}^{-1},\quad
		\partial^{R}_h=\bar{\mathcal{R}}\ast(\partial^{L}_h)\ast\bar{\mathcal{R}}^{-1}
	\end{equation}
	for
	$a, b\in A, \ h\in H$.
	Here,
	$(\partial^{L}_h)(a):=\langle h_{(1)},a\rangle h_{(2)},
	$
	$
	(\partial^{R}_h)(a):=h_{(1)}\langle h_{(2)},a\rangle$
are the left and right “differentiation operators”, regarded as maps $A \rightarrow H$ for any fixed $h$.

A quasitriangular Hopf algebra is a Hopf algebra $(H, m, \eta, \Delta, \epsilon, S)$
with a universal $R$-matrix $\R$,
an invertible element in $H\otimes H$ satisfying:
$$
\Delta^{\text{cop}}(a)=(\R)\Delta(a)(\R^{-1}),\,a \in H,\quad
(\Delta\otimes \text{id})(\R)=\R_{13}\R_{23},\quad
(\text{id}\otimes \Delta)(\R)=\R_{13}\R_{12},
$$
where
$\R_{12}=\R^{1}\otimes \R^{2}\otimes 1$,
$\R_{13}=\R^{1}\otimes 1\otimes \R^{2}$,
$\R_{23}=1\otimes\R^{1}\otimes \R^{2}$ for $\R=\R^{1}\otimes \R^{2}$.
An important example of quasitriangular Hopf algebras is
the (resp., $h$-adic) Drinfeld-Jimbo algebra $U_q(\mathfrak{g})$ (resp., $U_{h}(\mathfrak{g}))$ associated with a finite-dimensional complex simple Lie algebra $\mathfrak{g}$.
Denote the simple roots and the positive root system of $\mathfrak{g}$ by $\alpha_j$, $\Delta_{+}$, respectively.
Let $(a_{ij})$ be the Cartan matrix of $\mathfrak g$,
$(B_{ij})$ be the inverse of the matrix $\left(\frac{2}{(\alpha_j,\alpha_j)}a_{ij}\right)$.
The $q$-exponential function $\text{exp}_{q}x$ is defined by $\text{exp}_{q}x=\sum_{r=0}^{\infty}\frac{q^{r(r+1)/2}}{[r]_{q}!}x^{r}$
and $[n]_q=\frac{q^n-q^{-n}}{q-q^{-1}}$,
$[n]_q!=\prod_{k=1}^n [k]_q$.
Then an explicit form of the universal $R$-matrix of $U_{h}(\mathfrak{g})$ (\cite{klim}) is as follows:
\begin{equation}\label{imp2}
\R=\text{exp}\Bigl(h\sum\limits_{i,j}B_{ij}(H_{i}\otimes H_{j})\Bigr)\prod\limits_{\beta\in \Delta_{+}}\text{exp}_{q_{\beta}}\left((1-q_{\beta}^{-2})(E_{\beta}\otimes
F_{\beta})\right).
\end{equation}

Dually,
a coquasitriangular (or dual quasitriangular) Hopf algebra $(A,\mathfrak{r})$ is a Hopf algebra %
$A$ equipped with a convolution-invertible linear form $\mathfrak{r}: A\otimes
A \longrightarrow \mathbb{C}$ satisfying:
\begin{equation}\label{coquasi}
\mathfrak{r}\circ (\cdot\otimes \text{id})=\mathfrak{r}_{13}\ast\mathfrak{r}_{23},\quad
\mathfrak{r}\circ(\text{id}\otimes\cdot)=\mathfrak{r}_{13}\ast\mathfrak{r}_{12},\quad
\cdot\circ\tau=\mathfrak{r}\ast\cdot\ast\mathfrak{r}^{-1}.
\end{equation}
Here,
$\mathfrak{r}_{12}(a\otimes b\otimes c)=\mathfrak{r}(a\otimes b)\epsilon(c)$,
$\mathfrak{r}_{13}(a\otimes b\otimes c)=\mathfrak{r}(a\otimes c)\epsilon(b)$,
$\mathfrak{r}_{23}(a\otimes b\otimes c)=\mathfrak{r}(b\otimes c)\epsilon(a)$.

When the $R$-matrix $R$ satisfies the following \textit{FRT-condition} (cf. \cite{FRT}):
\begin{equation}\label{FRT}
(R^{-1})^{t_{1}}P(R^{t_{2}})^{-1}PK_{0}=\textrm{const}\cdot K_{0},
\end{equation}
there is a corresponding Hopf algebra $H_{R}$,
which is generated by  $1$ and $t^{i}_{j},\,\widetilde{t}^{i}_{j}, (i,j=1,\cdots,n)$,
with relations
\begin{gather*}
RT_{1}T_{2}=T_{2}T_{1}R,\quad
R^{t}\widetilde{T}_{1}\widetilde{T}_{2}=\widetilde{T}_{2}\widetilde{T}_{1}R^{t},\quad
(R^{t_{2}})^{-1}T_{1}\widetilde{T}_{2}=\widetilde{T}_{2}T_{1}(R^{t_{2}})^{-1};\\
\Delta(T)=T\otimes T,\quad
\Delta(\widetilde{T})=\widetilde{T}\otimes \widetilde{T},\quad
\epsilon(T)=\epsilon(\widetilde{T})=\textbf{I};\\
S(T)=(\tilde{T})^{t},\quad
S(\tilde{T})=DT^{t}D^{-1}.
\end{gather*}
Here, $t_1$ and $t_2$ denote the partial transposes in the first
and second tensor factors, respectively, while $t$ denotes the full
transpose.
$P$ denotes the matrix of the flip map.  When no confusion can arise, the same symbol denotes
both the flip and its matrix.
$(K_{0})^{ij}_{kl}:=\delta_{ij}\delta_{kl}$,
$D=\textrm{tr}_{2}(P((R^{t_{2}})^{-1})^{t_{1}})\in M_{n}(\mathbb{C})$, where the symbol $\operatorname {tr}_2$ denotes the partial trace
over the second tensor factor of $\mathbb C^n\otimes \mathbb C^n$.

For a finite-dimensional minuscule irreducible representation $T_V$ of $U_h(\mathfrak g)$ on
$V$, we proved in \cite {HH2} that the associated Hopf algebra $H_{R_V}$ is coquasitriangular,
where $R_V=(T_V\otimes T_V)(\mathscr R)$. Moreover, the FRT generators $m^{\pm
}{}^i_j:=S(l^{\pm }_{ij})$ associated with the $L$-functionals $l^{\pm}_{ij}$  of $T_V$ generate an extended
Hopf algebra $U_q^{\mathrm {ext}}(\mathfrak g)$ whose group-like elements are indexed by a
refinement of the weight lattice used for $U_q(\mathfrak g)$.
Furthermore, as a crucial ingredient in the generalized double-bosonization theory \cite{HH2}, we showed that
the pair $(U_q^{\mathrm {ext}}(\mathfrak g),H_{R_V})$ forms a weakly quasitriangular dual pair.

\subsubsection{Dually-paired braided groups}
Majid \cite {majid1} used the term \emph {braided groups} for bialgebras or Hopf algebras in a
braided tensor category, thereby distinguishing them from ordinary bialgebras and Hopf
algebras.
If a pair $(R,R’)$, called a \emph {Majid pair}, satisfies
\begin{gather}
	R_{12}R_{13}R^{\prime}_{23}=R^{\prime}_{23}R_{13}R_{12},\qquad
	R_{23}R_{13}R^{\prime}_{12}=R^{\prime}_{12}R_{13}R_{23},\label{*}\\
	(PR+\textbf{I})(PR'-\textbf{I})=0,\qquad
	R_{21}R^{\prime}_{12}=R^{\prime}_{21}R_{12},\label{**}
\end{gather}
then there is a braided group in the braided
category ${}^{H_R}\mathfrak{M}$,
named the \textit{braided vector algebra} $V(R',R)$,
generated by $1$, $\{e^{i}~|~i=1,\cdots,n\}$ with
relations
$$e^{i}e^{j}=\sum\limits_{a,b} R'{}^{ji}_{ab}e^{a}e^{b},\quad
\underline{\Delta}(e^{i})=e^{i}\otimes
1+1\otimes e^{i},\quad
\underline{\epsilon}(e^{i})=0,\quad
\underline{S}(e^{i})=-e^{i}.$$
Under the duality $\langle
f_{j},e^{i}\rangle=\delta_{ij}$,
the \textit{braided co-vector algebra} $V^{\vee}(R^{\prime},R_{21}^{-1})$ is a braided group in the braided
category $\mathfrak{M}^{H_R}$,
generated by $1$, $\{f_{j}\mid
j=1,\cdots,n\}$ with relations
$$f_{i}f_{j}=\sum\limits_{a,b}f_{b}f_{a}R'{}^{ab}_{ij},\quad
\underline{\Delta}(f_{i})=f_{i}\otimes
1+1\otimes f_{i},\quad
\underline{\epsilon}(f_{i})=0,\quad
\underline{S}(f_{i})=-f_{i}.$$
\begin{remark}\label{rem1}
In fact, the minimal polynomial equation $ \prod_{i}(PR-x_{i})=0
$ of the braiding $PR$ yields a Majid pair $(R, R')$  after $R$ is normalized at a
chosen eigenvalue $x_i$ of $PR$.
\end{remark}

\subsection{Generalized double-bosonization}
Starting from the dually-paired braided groups \(B^\star\) and B within the braided tensor category of $H$-modules over the quasitriangular Hopf algebra $H$, Majid constructed a Hopf-algebra structure on the tensor space \(B^\star\otimes H\otimes B\) by means of double-bosonization theory.
This Hopf algebra contains the two bosonizations \(B^\star\rtimes H\) and \(H\ltimes B\) as its Hopf subalgebras.
Based on the weakly quasitriangular dual pair $(U_q^{\textrm{ext}}(\mathfrak{g}), H_{R_{V}})$,
the above pair of $V^{\vee}(R^{\prime},R_{21}^{-1})$ and $V(R^{\prime},R)$ gives rise to the dually-paired braided groups \(B^\star,B\) in the braided category of \(U_q^{\textrm{ext}}(\mathfrak{g})\)-modules.
The generalized double-bosonization on the tensor space $V^{\vee}(R^{\prime},R_{21}^{-1})\otimes \widetilde{U_{q}^{ext}(\mathfrak g)}\otimes V(R^{\prime},R)$ is then obtained,
where
$\widetilde{U_{q}^{ext}(\mathfrak g)}:=U_{q}^{\textrm{ext}}(\mathfrak g)\otimes k[c,c^{-1}]$ denotes the central extension of $U_q^{\textrm{ext}}(\mathfrak{g})$.
Full details are given in \cite{HH2}.
\begin{theorem} $($\cite{HH2}$)$\label{cor1}
Let $R_{VV}$ be the $R$-matrix associated with a minuscule irreducible representation $T_{V}$ of $U_{q}(\mathfrak{g})$.
There exists a normalization constant $\lambda$ satisfying $\lambda R=R_{VV}$.
The new quantum group
$U=U(V^{\vee}(R^{\prime},R_{21}^{-1}),\widetilde{U_{q}^{ext}(\mathfrak g)},V(R^{\prime},R))$ obeys the following cross-commutation relations:
$$
cf_{i}=\lambda\, f_{i}\,c,\quad
e^{i}c=\lambda\, c\,e^{i},\quad
[c,m^{\pm}]=0,\quad
[e^{i},f_{j}]=\frac{(m^{+})^{i}_{j}\,c^{-1}-c\,(m^{-})^{i}_{j}}{q_{\ast}-q_{\ast}^{-1}},
$$
$$
e^{i}\,(m^{+})^{j}_{k}=R_{VV}{}^{ji}_{ab}\,(m^{+})^{a}_{k}\,e^{b},\quad
(m^{-})^{i}_{j}\,e^{k}=R_{VV}{}^{ki}_{ab}\,e^{a}\,(m^{-})^{b}_{j},$$
$$(m^{+})^{i}_{j}\,f_{k}=f_{b}\,(m^{+})^{i}_{a}\,R_{VV}{}^{ab}_{jk},\quad
f_{i}\,(m^{-})^{j}_{k}=(m^{-})^{j}_{b}\,f_{a}\,R_{VV}{}^{ab}_{ik},
$$
its counit is given by $\epsilon (e^{i})=\epsilon (f_{i})=0$,
and its comultiplication reads $$\Delta (c)=c\otimes c, \quad \Delta
(e^{i})=e^{a}\otimes (m^{+})^{i}_{a}\,c^{-1}+1\otimes e^{i}, \quad
\Delta (f_{i})=f_{i}\otimes 1+c\,(m^{-})^{a}_{i}\otimes f_{a}.$$
Here the factor \(q_*-q_*^{-1}\) is a normalization convention for the generators \(e^i\), introduced for computational convenience.
\end{theorem}

\section{Multi-Tensor-Product of Generalized Double-Bosonization}
Let $(H_{i},m_{H_i},1_{H_i},\Delta_{H_i},\epsilon_{H_i},S_{H_i})$ ($i=1,\cdots,n$)
be Hopf algebras.
The permutation map $P$ on $V\otimes V$ for a vector space $V$ extends naturally to any tensor power $V^{\otimes n}$ ($n\geq 2$).
Then each operator $P_{ij}$ is the image of the corresponding transposition $(ij)\in \mathfrak{G}_{n}$ under the natural action of the symmetric group $\mathfrak{G}_{n}$
on $V^{\otimes n}$ by permuting the tensor factors.
Consider the following two elements of $\mathfrak{G}_{2n}$:
\begin{gather*}
\sigma=\left(
\begin{array}{cccccccccc}
1&2&3&4&\cdots&2k-1&2k&\cdots&2n-1&2n\\
1&n+1&2&n+2&\cdots&k&n+k&\cdots&n&2n
\end{array}
\right),
\\
\tau=\left(
\begin{array}{cccccccccccc}
1&2&\cdots&k&\cdots&n&n+1&n+2&\cdots&n+k&\cdots&2n\\
1&3&\cdots&2k-1&\cdots&2n-1&2&4&\cdots&2k&\cdots&2n
\end{array}
\right).
\end{gather*}
Then $\sigma\tau=\tau\sigma=1$.
The associated operators acting on $V^{\otimes 2n}$ are denoted by $\mathcal{X}_{nn}, \mathcal{X}_{nn}^{-1}$,
where
\begin{gather*}
\mathcal{X}_{nn}=(P_{n,n+1}P_{n+1,n+2}\cdots P_{2n-2,2n-1})(P_{n-1,n}P_{n,n+1}\cdots P_{2n-4,2n-3})\cdots(P_{34}P_{45})(P_{23}),\\
\mathcal{X}_{nn}^{-1}=(P_{2n-2,2n-1})(P_{2n-4,2n-3}P_{2n-3,2n-2})\cdots\cdots(P_{45}\cdots P_{n,n+1}P_{n+1,n+2})(P_{23}\cdots P_{n-1,n}P_{n,n+1}).
\end{gather*}
The corresponding structure maps of $H_{1}\otimes \cdots\otimes H_{n}$ are
\begin{gather}\label{hopf}
m_{H_{1}\otimes \cdots\otimes H_{n}}=(m_{H_1}\otimes\cdots\otimes m_{H_n})\circ\mathcal{X}_{nn}^{-1},\quad
1_{H_{1}\otimes \cdots\otimes H_{n}}=1_{H_1}\otimes\cdots\otimes 1_{H_{n}};
\\
\Delta_{H_{1}\otimes \cdots\otimes H_{n}}=\mathcal{X}_{nn}\circ(\Delta_{H_1}\otimes\cdots\otimes\Delta_{H_n}),\quad
\epsilon_{H_{1}\otimes \cdots\otimes H_{n}}=\epsilon_{H_{1}}\otimes \cdots\otimes\epsilon_{H_{n}};
\\
S_{H_{1}\otimes \cdots\otimes H_{n}}=S_{H_1}\otimes\cdots\otimes S_{H_n}.
\end{gather}
The Hopf algebra on $H_1\otimes\cdots\otimes H_n$ is equipped with its tensor-product structure.

\subsection{Multi-tensor-product of (co-)quasitriangular Hopf algebras}
If, for each $i=1,\cdots$, $n$, the Hopf algebra $H_i$ is quasitriangular with universal
$R$-matrix $\mathscr R_{H_i}$, then the following lemma holds.
\begin{lemma}\label{lemqua}
The multi-tensor-product Hopf algebra $H_{1}\otimes\cdots\otimes H_n$ is also quasitriangular,
whose universal $R$-matrix is given by
$\R_{H_{1}\otimes\cdots\otimes H_n}=\mathcal{X}_{nn}\circ(\R_{H_1}\otimes\cdots\otimes\R_{H_n})$.
\end{lemma}
Similarly,
suppose \((A_i,\mathfrak{r}_{A_i})\) is a coquasitriangular Hopf algebra for each \(i=1,\dots,n\).
Denote by $A_1\otimes\cdots\otimes A_n$ the multi-tensor-product Hopf algebra.
 The following lemma is the dual counterpart of Lemma \ref{lemqua}. 
\begin{lemma}\label{lemco} 
The Hopf algebra $A_1\otimes\cdots\otimes A_n$ is also coquasitriangular, with the convolution-invertible linear form
$\mathfrak{r}_{A_1\otimes\cdots\otimes A_n}
=(\mathfrak{r}_{A_1}\otimes\cdots\otimes\mathfrak{r}_{A_n})\circ\mathcal{X}_{nn}^{-1}$.
\end{lemma}

Let $T_i$ (resp. $\rho_i)$ be a (co-)representation of the (co-)quasitriangular Hopf algebra $H_{i}$ (resp. $A_i)$ on an $m_i$-dimensional vector space $V_i$,
$i=1,\cdots,n$.
It is straightforward to verify that
$(T_1\otimes\cdots\otimes T_n)\circ \mathcal{X}_{nn}^{-1}$ is a representation of the quasitriangular Hopf algebra
$H_1\otimes\cdots\otimes H_n$ on the $(m_1\cdots m_n)$-dimensional vector space $V_1\otimes\cdots \otimes V_n$,
which we denote by $T_{V_1\otimes\cdots \otimes V_n}$. Dually,
$\mathcal{X}_{nn}\circ(\rho_1\otimes\cdots\otimes \rho_n)$ is a corepresentation of the coquasitriangular Hopf algebra
$(A_1\otimes\cdots\otimes A_n)$ on $V_1\otimes\cdots \otimes V_n$,
denoted by $\rho_{_{V_1\otimes\cdots \otimes V_n}}$.

Owing to the quasitriangularity of each \(H_i\), fixing a basis \(S_i\) of \(V_i\) yields an invertible $R$-matrix
\(R_{V_i}=(T_i\otimes T_i)(\mathscr{R}_{H_i})\).
With respect to the basis \(S_1\otimes\cdots\otimes S_n\) of the \(H_1\otimes\cdots\otimes H_n\)-module \(V_1\otimes\cdots \otimes V_n\), we similarly obtain an invertible \((m_1\cdots m_n)^2\times (m_1\cdots m_n)^2\) $R$-matrix from the universal $R$-matrix \(\mathscr{R}_{H_1\otimes\cdots \otimes H_n}\) given in Lemma \ref{lemqua}; this matrix is denoted by \(R_{V_1\otimes\cdots \otimes V_n}\).

Let $(R_{V_1\otimes\cdots \otimes V_n})^{(i_1i_{2}\cdots i_n)(j_1j_{2}\cdots j_n)}_{(k_1k_{2}\cdots k_n)(l_1l_{2}\cdots l_n)}$ denote the entry located at row
$(i_1i_{2}\cdots i_n)(j_1j_{2}\cdots j_n)$ and column $(k_1k_{2}\cdots k_n)(l_1l_{2}\cdots l_n)$.
For notational simplicity, we write \(\underline{i}\) for the $n$-tuple $(i_{1},\cdots,i_{n})$,
where $i_{k}\in \mathbb{Z_{+}}$, for $k=1, \cdots, n$.
The above matrix entry may be abbreviated as
$(R_{V_1\otimes\cdots \otimes V_n})^{\underline{i}\,\underline{j}}_{\underline{k}\,\underline{l}}$.
Direct computation of each entry of \(R_{V_1\otimes\cdots\otimes V_n}\) is rather involved, whereas Proposition \ref{rmatrix} below provides a compact description.
\begin{proposition}\label{rmatrix}
The matrix entries satisfy
\((R_{V_1\otimes\cdots \otimes V_n})^{\underline{i}\,\underline{j}}_{\underline{k}\,\underline{l}}
=(R_{V_1})^{i_1j_1}_{k_1l_1}(R_{V_2})^{i_2j_2}_{k_2l_2}\cdots(R_{V_n})^{i_nj_n}_{k_nl_n}.\)
\end{proposition}
\begin{proof}
In terms of the universal $R$-matrix $\R_{V_1\otimes\cdots V_n}$ from Lemma \ref{lemqua},
we have that
\begin{equation}\label{r1}
\left.
\begin{array}{rl}
(T_{V_1\otimes\cdots \otimes V_n}\otimes& T_{V_1\otimes\cdots \otimes V_n})(\R_{H_1\otimes\cdots \otimes H_n})\Bigl((v_{k_1}\otimes v_{k_2}\cdots\otimes v_{k_n})\otimes (v_{l_1}\otimes
v_{l_2}\otimes\cdots\otimes v_{l_n})\Bigr)\\
=&(R_{V_1\otimes\cdots \otimes V_n})^{(i_1i_{2}\cdots i_n)(j_1j_{2}\cdots j_n)}_{(k_1k_{2}\cdots k_n)(l_1l_{2}\cdots l_n)}\Bigl((v_{i_1}\otimes v_{i_2}\otimes\cdots\otimes
v_{i_n})\otimes (v_{j_1}\otimes v_{j_2}\otimes\cdots\otimes v_{j_n})\Bigr)\\
=&(R_{V_1\otimes\cdots \otimes V_n})^{\underline{i}\,\underline{j}}_{\underline{k}\,\underline{l}}\Bigl((v_{i_1}\otimes v_{i_2}\otimes\cdots\otimes v_{i_n})\otimes (v_{j_1}\otimes
v_{j_2}\otimes\cdots\otimes v_{j_n})\Bigr).
\end{array}
\right.
\end{equation}
On the other hand,
\begin{equation}\label{r2}
\left.
\begin{array}{rl}
(T_{V_1\otimes\cdots \otimes V_n}&\otimes T_{V_1\otimes\cdots \otimes V_n})(\R_{H_1\otimes\cdots H_n})\Bigl((v_{k_1}\otimes\cdots\otimes v_{k_n})\otimes (v_{l_1}\otimes\cdots\otimes v_{l_n})\Bigr)\\
&=(T_{V_1\otimes\cdots \otimes V_n})(\R_{H_1}^{1}\otimes\cdots\otimes\R_{H_n}^{1})(v_{k_1}\otimes\cdots\otimes v_{k_n})\,\otimes\\
&\qquad\qquad\quad\;(T_{V_1\otimes\cdots \otimes V_n})(\R_{H_1}^{2}\otimes\cdots\otimes\R_{H_n}^{2})(v_{l_1}\otimes\cdots\otimes v_{l_n})\\
&=T_{V_1}(\R_{H_1}^{1})(v_{k_1})\otimes\cdots\otimes T_{V_n}(\R_{H_n}^{1})(v_{k_n})\otimes
T_{V_1}(\R_{H_1}^{2})(v_{l_1})\otimes\cdots\otimes T_{V_n}(\R_{H_n}^{2})(v_{l_n})\\
&=\mathcal X_{nn}\circ T_{V_1}(\R_{H_1}^{1})(v_{k_1})\otimes T_{V_1}(\R_{H_1}^{2})(v_{l_1})\otimes \cdots\otimes T_{V_n}(\R_{H_n}^{1})(v_{k_n})\otimes T_{V_n}(\R_{H_n}^{2})(v_{l_n})\\
&=\mathcal X_{nn}\circ(T_1\otimes T_1)(\R_{H_1})(v_{k_1}\otimes v_{l_1})\otimes\cdots\otimes(T_n\otimes T_n)(\R_{H_n})(v_{k_n}\otimes v_{l_n})
\\
&=\mathcal X_{nn}\circ(R_{V_1})^{i_1j_1}_{k_1l_1}(v_{i_1}\otimes v_{j_1})\otimes\cdots\otimes
(R_{V_n})^{i_nj_n}_{k_nl_n}(v_{i_n}\otimes v_{j_n})\\
&=\mathcal X_{nn}\circ(R_{V_1})^{i_1j_1}_{k_1l_1}\cdots(R_{V_n})^{i_nj_n}_{k_nl_n}\Bigl((v_{i_1}\otimes v_{j_1})\otimes\cdots\otimes(v_{i_n}\otimes v_{j_n})\Bigr)\\
&=(R_{V_1})^{i_1j_1}_{k_1l_1}(R_{V_2})^{i_2j_2}_{k_2l_2}\cdots(R_{V_n})^{i_nj_n}_{k_nl_n}\mathcal X_{nn}\circ\Bigl((v_{i_1}\otimes v_{j_1})\otimes\cdots\otimes(v_{i_n}\otimes v_{j_n})\Bigr)\\
&=(R_{V_1})^{i_1j_1}_{k_1l_1}(R_{V_2})^{i_2j_2}_{k_2l_2}\cdots(R_{V_n})^{i_nj_n}_{k_nl_n}\Bigl((v_{i_1}\otimes \cdots\otimes v_{i_n})\otimes (v_{j_1}\otimes \cdots\otimes v_{j_n})\Bigr).
\end{array}
\right.
\end{equation}
By virtue of \eqref{r1} and \eqref{r2},
we obtain $(R_{V_1\otimes\cdots \otimes V_n})^{\underline{i}\,\underline{j}}_{\underline{k}\,\underline{l}}
=(R_{V_1})^{i_1j_1}_{k_1l_1}(R_{V_2})^{i_2j_2}_{k_2l_2}\cdots(R_{V_n})^{i_nj_n}_{k_nl_n}$.
\end{proof}

The matrix $PR_{V_1\otimes\cdots\otimes V_n}$ represents the operator
$\mathcal{P}\circ(T_{V_1\otimes\cdots\otimes V_n}\otimes T_{V_1\otimes\cdots\otimes V_n})(\R_{H_1\otimes\cdots\otimes H_n})$ acting on tensor space
$(V_1\otimes\cdots\otimes V_n)\otimes(V_1\otimes\cdots\otimes V_n)$.
Here,
the permutation operator $\mathcal{P}$ corresponds to the following element of the symmetric group $\mathfrak{G}_{2n}$:
\begin{gather*}
\left(
\begin{array}{ccccccccccc}
1&2&\cdots&k&\cdots&n&n+1&\cdots&n+k&\cdots&2n\\
n+1&n+2&\cdots&n+k&\cdots&2n&1&\cdots&k&\cdots&n
\end{array}
\right).
\end{gather*}
Starting from any invertible $R$-matrix $R$,
we need to find the Majid pair $(R, R')$ satisfying conditions \eqref{*} and \eqref{**}.
We therefore determine the characteristic polynomial of $PR_{V_1\otimes\cdots\otimes V_n}$.
For each $i$,
the braiding associated with $R_{V_i}$ is the operator $PR_{V_i}$ on $V_i\otimes V_i$, where $P$ is the flip operator.
Let $I_{R_{V_i}}$ be the set of its eigenvalues; equivalently, $\lambda\in k$ satisfies
$\det(PR_{V_i}-\lambda\textrm{I})=0$ if and only if $\lambda\in I_{R_{V_i}}$.
Then we have
\begin{proposition}\label{minimal}
The eigenvalues of the braiding $PR_{V_1\otimes\cdots\otimes V_n}$ are of the form $\lambda_{i_1}\cdots\lambda_{i_n}$, where $\lambda_{i_j}\in I_{R_{V_j}}$ for every $j\in\{1,\cdots,n\}$.
Furthermore, if each $PR_{V_i}$ is diagonalizable, then so is $PR_{V_1\otimes\cdots V_n}$.
\end{proposition}
\begin{proof}
Let $\mathcal{P}_{i,n+i}$ be the operator on $(V_1\otimes\cdots\otimes V_n)\otimes(V_1\otimes\cdots\otimes V_n)$ corresponding to the element
\begin{gather*}
\left(
\begin{array}{cccccccccccccc}
1&\cdots&i-1&i&i+1&\cdots&n&n+1&\cdots&n+i-1&n+i&n+i+1&\cdots&2n\\
1&\cdots&i-1&n+i&i+1&\cdots&n&n+1&\cdots&n+i-1&i&n+i+1&\cdots&2n
\end{array}
\right)
\end{gather*}
of $\mathfrak{G}_{2n}$ for any $i\in\{1,2,\cdots,n\}$.
Obviously,
$\mathcal{P}_{i,n+i}$ and $\mathcal{P}_{j,n+j}$ commute if the indices $i, j$ are distinct.
Accordingly,
let \(T_{i,n+i}\) denote the operator given below
\begin{equation*}
\left.
\begin{array}{rl}
\mathcal{P}_{i,n+i}\circ[(\text{id}_{V_1\otimes\cdots \otimes V_{i-1}}\otimes T_{V_i}(\R_{H_i}^1)\otimes\text{id}_{V_{i+1}\otimes\cdots\otimes
V_n})\otimes(\text{id}_{V_1\otimes\cdots \otimes V_{i-1}}\otimes T_{V_i}(\R_{H_i}^2)\otimes\text{id}_{V_{i+1}\otimes\cdots\otimes V_n})]
\end{array}
\right.
\end{equation*}
acting on vector space $(V_1\otimes\cdots\otimes V_n)\otimes(V_1\otimes\cdots\otimes V_n)$.
Then we have
\begin{equation}\label{r3}
\mathcal{P}\circ(T_{V_1\otimes\cdots\otimes V_n}\otimes T_{V_1\otimes\cdots\otimes V_n})(\R_{H_1\otimes\cdots\otimes H_n})
=T_{1,n+1}\circ\cdots\circ T_{i,n+i}\circ\cdots\circ T_{n,n+n}.
\end{equation}
It is easy to see that operators $T_{i,n+i}$ and $T_{j,n+j}$ also commute when the indices $i, j$ are distinct.
Furthermore,
as operators acting on the vector space $(V_1\otimes\cdots\otimes V_n)\otimes(V_1\otimes\cdots\otimes V_n)$,
we know that
$T_{i,n+i}={\mathcal{P}}_{i+1,n+i}\circ\Bigl(\text{id}\otimes \left((P\circ (T_{V_i}\otimes T_{V_i}))(\R_{H_i})\right)\otimes\text{id}\Bigr)\circ\mathcal{P}_{i+1,n+i}$.
Since the operator ${\mathcal{P}}_{i+1,n+i}$ is invertible,
and ${\mathcal{P}}_{i+1,n+i}\circ{\mathcal{P}}_{i+1,n+i}=\textrm{id}$,
the matrix \(M_{i,n+i}\) associated with \(T_{i,n+i}\) is conjugate to
\(\mathrm{I}^{\otimes(i-1)}\otimes PR_{V_i}\otimes \mathrm{I}^{\otimes (2n-i-1)}\).
Consequently, \(M_{i,n+i}\) shares the same minimal polynomial as \(PR_{V_i}\).
Combining this identity with \eqref{r3}, we obtain
\(PR_{V_1\otimes\cdots\otimes V_n}=M_{1,n+1}M_{2,n+2}\cdots M_{i,n+i}\cdots M_{n,2n},\)
where \(M_{i,n+i}\) commutes with \(M_{j,n+j}\) whenever \(i\neq j\).
The Proposition now follows from standard facts in linear algebra.
\end{proof}

\subsection{Weakly quasitriangular dual pairs}
Given dual pairs of bialgebras or Hopf algebras \((H_i,A_i)\) equipped with the pairing \(\langle\,,\,\rangle_i\),
the bialgebra (or Hopf-algebra) pair \(H_1\otimes\cdots\otimes H_n\) and \(A_1\otimes\cdots\otimes A_n\) likewise forms a dual pairing under
\(\langle h_1\otimes\cdots\otimes h_n,\,a_1\otimes\cdots\otimes a_n\rangle:=\langle h_1,a_1\rangle_1\cdots\langle h_n,a_n\rangle_n.\)
Correspondingly, we also have a tensor-product counterpart among weakly quasitriangular dual pairs.
Suppose we have a family of weakly quasitriangular dual pairs \((H_i,A_i)\) (\(i=1,\dots,n\)) of bialgebras (or Hopf algebras).
Let \(\mathcal{R}_i\), \(\bar{\mathcal{R}}_i\) denote their respective convolution-invertible algebra/anti-coalgebra maps,
and, for each fixed \(h_i\in H_i\), let \(\partial^L h_i\), \(\partial^R h_i\) stand for the left- and right-hand “differentiation operators”.
The we have
\begin{lemma}\label{weak}
The bialgebra- or Hopf-algebra pair \(H_1\otimes\cdots\otimes H_n\) and \(A_1\otimes\cdots\otimes A_n\) forms a weakly quasitriangular dual pair,
with convolution-invertible algebra/anti-coalgebra maps
\(\mathcal{R}=\mathcal{R}_{1}\otimes\cdots\otimes\mathcal{R}_{n},\qquad
\bar{\mathcal{R}}=\bar{\mathcal{R}}_{1}\otimes\cdots\otimes\bar{\mathcal{R}}_{n}\colon\quad
A_1\otimes\cdots\otimes A_n \longrightarrow H_1\otimes\cdots\otimes H_n,\)
which obey
\begin{gather}
\langle\bar{\mathcal{R}}(a_1\otimes\cdots\otimes a_n),b_1\otimes\cdots\otimes b_n\rangle
=\langle\mathcal{R}^{-1}(b_1\otimes\cdots\otimes b_n),a_1\otimes\cdots\otimes a_n\rangle,\\
\partial^{R}_{h_1\otimes\cdots\otimes h_n}=\mathcal{R} \ast(\partial^{L}_{h_1\otimes\cdots\otimes h_n})\ast\mathcal{R}^{-1},\quad
\partial^{R}_{h_1\otimes\cdots\otimes h_n}=\bar{\mathcal{R}}\ast(\partial^{L}_{h_1\otimes\cdots\otimes h_n})\ast\bar{\mathcal{R}}^{-1}.
\end{gather}
Here, for any fixed \(h_1\otimes\cdots\otimes h_n\), the left- and right-hand “differentiation operators”,
viewed as maps \(A_1\otimes\cdots\otimes A_n \longrightarrow H_1\otimes\cdots\otimes H_n\), satisfy
$\partial^{R}_{h_1\otimes\cdots\otimes h_n}=\partial^{R}_{h_1}\otimes\cdots\otimes \partial^{R}_{h_n}$,
$\partial^{L}_{h_1\otimes\cdots\otimes h_n}=\partial^{L}_{h_1}\otimes\cdots\otimes \partial^{L}_{h_n}$.
\end{lemma}
We now restrict our attention to quasitriangular Hopf algebras \(U_q(\mathfrak{g})\) (or \(U_h(\mathfrak{g})\)) and their finite‑dimensional minuscule irreducible representations.
This standing assumption will be in force unless otherwise specified.
Let $\{\mathfrak{g}_1, \cdots, \mathfrak{g}_n\}$ be a set of complex simple Lie algebras.
We obtain a universal $R$-matrix of $U_h(\mathfrak{g}_1\oplus\cdots\oplus\mathfrak{g}_n)$ in view of the isomorphism $U_h(\mathfrak{g}_1\oplus\cdots\oplus\mathfrak{g}_n)\simeq U_h(\mathfrak{g}_1)\otimes\cdots\otimes U_h(\mathfrak{g}_n)$,
denoted still by $\R$.
Associated with irreducible $U_{h}(\mathfrak{g}_i)$-representation $T_{V_i}$ on $m_i$-dimensional vector space $V_i$,
and a basis $\{x^{(i)}_j, j=1,\cdots,m_i\}$ of $V_i$,
we get the corresponding $L$-functionals $L^{\pm}$ being subordinate to $U_{h}(\mathfrak{g}_i)$,
that is,
$$
(\textrm{id}\otimes T_{V_i})\bigl(\R^{1}_{U_h(\mathfrak{g_i})}\otimes\R^{2}_{U_h(\mathfrak{g_i})}\bigr)(1\otimes x^{(i)}_{j_i})
=\R^{1}_{U_h(\mathfrak{g_i})}\otimes T_{V_i}(\R^{2}_{U_h(\mathfrak{g_i})})(x^{(i)}_{j_i})
=\sum_{k_i=1}^{m_i}l^{+}_{k_i,j_i}\otimes x^{(i)}_{k_i}.
$$
Here $k_i$ is the basis index running from $1$ to $m_i$, and the displayed expression is summed over $k_i$.
Denote by $x_{\underline{j}}$ the vector $x^{(1)}_{j_1}\otimes\cdots\otimes x^{(n)}_{j_n}$ of $V_1\otimes\cdots\otimes V_n$, and, as $\underline{j}$ varies, these vectors form a basis of $V_1\otimes\cdots\otimes V_n$.
Then
$$(\textrm{id}\otimes T_{V_1\otimes\cdots\otimes V_n})(\R)(1\otimes x_{\underline{j}})=\sum\limits_{\underline{i}}l^{+}_{\underline{i}\underline{j}}\otimes x_{\underline{i}},\quad
(T_{V_1\otimes\cdots\otimes V_n}\otimes \textrm{id})(\R^{-1})(x_{\underline{j}}\otimes 1)=\sum_{\underline{i}}x_{\underline{i}}\otimes l^{-}_{\underline{i}\underline{j}}.$$
These $l^{\pm}_{\underline{i}\underline{j}}\in U_{h}(\mathfrak{g}_1\oplus\cdots\oplus\mathfrak{g}_n)$ are $L$-functionals associated with the representation $T_{V_1\otimes\cdots\otimes V_n}$.
With these in hand, we have the following equalities up to an isomorphism:
\begin{equation*}
\begin{split}
&(\textrm{id}\otimes T_{V_1\otimes\cdots\otimes V_n})(\R)(1\otimes x_{\underline{j}})\\
&\ =(\textrm{id}\otimes T_{V_1{\otimes}\cdots{\otimes} V_n})((\R^1_{U_h(\mathfrak{g}_1)}{\otimes}\cdots{\otimes}\R^1_{U_h(\mathfrak{g}_n)})
\otimes(\R^2_{U_h(\mathfrak{g}_1)}{\otimes}\cdots{\otimes}\R^{2}_{U_h(\mathfrak{g}_n)}))(1\otimes x^{(1)}_{j_1}\otimes\cdots\otimes x^{(n)}_{j_n})\\
&\ =\R^{1}_{U_h(\mathfrak{g}_1)}\otimes\cdots\otimes\R^{1}_{U_h(\mathfrak{g}_n)}\otimes T_{V_1}(\R^{2}_{U_h(\mathfrak{g}_1)})(x^{(1)}_{j_1})
\otimes\cdots\otimes T_{V_n}(\R^{2}_{U_h(\mathfrak{g}_n)})(x^{(n)}_{j_n})\\
&\ =\sum_{k_1=1}^{m_1}\cdots\sum_{k_n=1}^{m_n}l^{+}_{k_1,j_1}\otimes\cdots\otimes l^{+}_{k_n,j_n}\otimes x^{(1)}_{k_1}\otimes\cdots\otimes x^{(n)}_{k_n}\\
&\ =\sum_{k_1=1}^{m_1}\cdots\sum_{k_n=1}^{m_n}l^{+}_{k_1,j_1}\otimes\cdots\otimes l^{+}_{k_n,j_n}\otimes x_{\underline{k}}.
\end{split}
\end{equation*}
Thus we have $l^{+}_{\underline{k}\underline{j}}=l^{+}_{k_1,j_1}\otimes\cdots\otimes l^{+}_{k_n,j_n}$,
and the relation $l^{-}_{\underline{k}\underline{j}}=l^{-}_{k_1,j_1}\otimes\cdots\otimes l^{-}_{k_n,j_n}$ can be obtained similarly.
Combining this with the relations between the FRT-generators $m^{\pm}$ and $L$-functionals,
we prove the following Proposition \ref{pm}.
\begin{proposition}\label{pm}
In terms of FRT-generators $m^{\pm}$ associated with $U^{\textrm{ext}}_q(\mathfrak{g_i})$-representation $T_{V_i}$,
we obtain (up to isomorphism) FRT-generators $(m^{\pm})^{\underline{k}}_{\underline{j}}$ of $U_q^{\textrm{ext}}(\mathfrak{g}_1\oplus\cdots\oplus\mathfrak{g}_n)$ corresponding to its representation  $T_{V_1\otimes\cdots\otimes V_n}$ as
$$(m^{+})^{\underline{k}}_{\underline{j}}=(m^{+})^{k_1}_{j_1}\otimes\cdots\otimes (m^{+})^{k_n}_{j_n}
\quad\mbox{\textrm{and}} \quad(m^{-})^{\underline{k}}_{\underline{j}}=(m^{-})^{k_1}_{j_1}\otimes\cdots\otimes (m^{-})^{k_n}_{j_n}.$$
\end{proposition}
There is a weakly quasitriangular dual pair $(U_q^{\textrm{ext}}(\mathfrak{g}_i),H_{R_{V_i}})$ 
for some minuscule representation $V_i$ of $U_q(\mathfrak{g}_i)$, $i=1,\cdots,n$ in \cite{HH2},
then we have a weakly quasitriangular dual pair $(U_q^{\textrm{ext}}(\mathfrak{g}_1)\otimes\cdots\otimes U_q^{\textrm{ext}}(\mathfrak{g}_n),H_{R_{V_1}}\otimes\cdots\otimes H_{R_{V_n}})$ by Lemma \ref{weak},
and a dual pair of braided co-vector algebra $V(R_i',R_i)$ and  braided vector algebra $V(R_i',(R_i)^{-1}_{21})$ in the braided $H_{R_{V_i}}$-comodules category, respectively.

However, these objects live in different categories of braided \(H_{R_{V_i}}\)-comodules (resp., \(U_q^{\text{ext}}(\mathfrak{g}_i)\)-modules).
Accordingly, the coquasitriangular Hopf algebra \(H_{R_{V_1}}\otimes\cdots\otimes H_{R_{V_n}}\), together with the braided-group pair
\(V(R_1',(R_1)_{21}^{-1})\otimes\cdots\otimes V(R_n',(R_n)_{21}^{-1})\quad\text{and}\quad V(R_1',R_1)\otimes\cdots\otimes V(R_n',R_n),\)
does not yield a suitable candidate for our purposes.
It is thus necessary to find a replacement for \(H_{R_{V_1}}\otimes\cdots\otimes H_{R_{V_n}}\) such that it, jointly with \(U_q^{\text{ext}}(\mathfrak{g}_1)\otimes\cdots\otimes U_q^{\text{ext}}(\mathfrak{g}_n)\), forms another weakly quasitriangular dual pair.

To this end, we return to the relevant $R$-matrices.
Given a finite-dimensional minuscule irreducible representation \(T_{V_i}\) of \(U_{q}(\mathfrak{g}_i)\) on the \(m_i\)-dimensional vector space \(V_i\),
we have the associated $R$-matrix \(R_{V_i}\), for which \(R_{V_i}^{t_{2}}\) is invertible, and there exists a Majid pair \((R_i,R_i')\) satisfying conditions \eqref{*}, \eqref{**}, and so on.
Correspondingly, the coquasitriangular Hopf algebra \(H_{R_{V_i}}\) is well-defined.
We now consider the $R$-matrix \(R_{V_1\otimes\cdots\otimes V_n}\).
Whenever each \(R_{V_i}^{t_2}\) is invertible, we define the matrix given in \eqref{N} as follows:
\begin{equation}\label{N}
\mathcal{N}^{\underline{i}\,\underline{j}}_{\underline{k}\,\underline{l}}=\mathcal{N}^{(i_1i_2\cdots i_n)(j_1j_2\cdots j_n)}
_{(k_1k_2\cdots k_n)(l_1l_2\cdots l_n)}=((R_{V_1}^{t_2})^{-1})^{i_1j_1}_{k_1l_1}((R_{V_2}^{t_2})^{-1})^{i_2j_2}_{k_2l_2}
\cdots((R_{V_n}^{t_2})^{-1})^{i_nj_n}_{k_nl_n},
\end{equation}
and we then have
\begin{equation}\label{t2}
\left.
\begin{array}{rl}
&\sum\limits_{\underline{a},\underline{b}}((R_{V_1\otimes\cdots\otimes V_n})^{t_2})^{\underline{i}\,\underline{j}}_{\underline{a}\,\underline{b}}
\mathcal{N}^{\underline{a},\underline{b}}_{\underline{k}\,\underline{l}}\\
&\quad=\sum\limits_{\underline{a},\underline{b}}(R_{V_1}^{t_2})^{i_1j_1}_{a_1b_1}(R_{V_2}^{t_2})^{i_2j_2}_{a_2b_2}
\cdots(R_{V_n}^{t_2})^{i_nj_n}_{a_nb_n}((R_{V_1}^{t_2})^{-1})^{a_1b_1}_{k_1l_1}((R_{V_2}^{t_2})^{-1})^{a_2b_2}_{k_2l_2}
\cdots((R_{V_n}^{t_2})^{-1})^{a_nb_n}_{k_nl_n}\\
&\quad=\sum\limits_{\underline{a},\underline{b}}(R_{V_1}^{t_2})^{i_1j_1}_{a_1b_1}((R_{V_1}^{t_2})^{-1})^{a_1b_1}_{k_1l_1}
(R_{V_2}^{t_2})^{i_2j_2}_{a_2b_2}((R_{V_2}^{t_2})^{-1})^{a_2b_2}_{k_2l_2}
\cdots(R_{V_n}^{t_2})^{i_nj_n}_{a_nb_n}((R_{V_n}^{t_2})^{-1})^{a_nb_n}_{k_nl_n}\\
&\quad=\delta_{i_1k_1}\delta_{j_1l_1}\delta_{i_2k_2}\delta_{j_2l_2}\cdots\delta_{i_nk_n}\delta_{j_nl_n}
=\delta_{\underline{i}\,\underline{k}}\delta_{\underline{j}\,\underline{l}}.
\end{array}
\right.
\end{equation}
The companion identity $\sum\limits_{\underline{a}\,\underline{b}}\mathcal{N}^{\underline{i}\,\underline{j}}_{\underline{a}\,\underline{b}}
((R_{V_1\otimes\cdots\otimes
V_n})^{t_2})^{\underline{a}\,\underline{b}}_{\underline{k}\,\underline{l}}=\delta_{\underline{i}\,\underline{k}}\delta_{\underline{j}\,\underline{l}}$
follows by an analogous argument.
We therefore conclude that the matrix \(R_{V_1\otimes\cdots\otimes V_n}^{t_2}\) is invertible.

On the other hand,
if each $R_{V_i}$ satisfies {\it FRT-condition} \eqref{FRT} described in Subsection \ref{gel},
then we have
\begin{equation}\label{hopf1}
\left.
\begin{array}{rl}
&\left[(R_{V_1\otimes\cdots\otimes V_n}^{-1})^{t_1}P(R_{V_1\otimes\cdots\otimes V_n}^{t_2})^{-1}PK_0\right]^{\underline{i}\,\underline{j}}_{\underline{k}\,\underline{l}}\\
&\quad=\sum\limits_{
\underline{a},\underline{b},\underline{c},\underline{d}}
((R_{V_1\otimes\cdots\otimes V_n}^{-1})^{t_1})^{\underline{i}\,\underline{j}}_{\underline{a}\,\underline{b}}
(P(R_{V_1\otimes\cdots\otimes
V_n}^{t_2})^{-1}P)^{\underline{a}\,\underline{b}}_{\underline{c}\,\underline{d}}(K_0)^{\underline{c}\,\underline{d}}_{\underline{k}\,\underline{l}}\\
&\quad=\sum\limits_{\underline{a},\underline{b},\underline{c},\underline{d}}
[((R_{V_1}^{{-}1})^{t_1})^{i_1j_1}_{a_1b_1}(P(R_{V_1}^{t_2})^{{-}1}P)^{a_1b_1}_{c_1d_1}(K_0)^{c_1d_1}_{k_1l_1}]
{\cdots}[((R_{V_n}^{{-}1})^{t_1})^{i_nj_n}_{a_nb_n}(P(R_{V_n}^{t_2})^{{-}1}P)^{a_nb_n}_{c_nd_n}(K_0)^{c_nd_n}_{k_nl_n}]\\
&\quad=[(R_{V_1}^{-1})^{t_1}P(R_{V_1}^{t_2})^{-1}PK_0]^{i_1j_1}_{k_1l_1}\cdots[(R_{V_n}^{-1})^{t_1}P(R_{V_n}^{t_2})^{-1}PK_0])^{i_nj_n}_{k_nl_n}\\
&\quad=\text{const}\cdot\delta_{i_1j_1}\delta_{k_1l_1}\cdots\delta_{i_nj_n}\delta_{k_nl_n}
=\text{const}\cdot\delta_{\underline{i}\,\underline{j}}\delta_{\underline{k}\,\underline{l}},
\end{array}
\right.
\end{equation}
which proves that the {\it FRT condition}: $(R_{V_1\otimes\cdots\otimes V_n}^{-1})^{t_1}P(R_{V_1\otimes\cdots\otimes V_n}^{t_2})^{-1}PK_0=\text{const}\cdot K_0$ holds.
This yields the first assertion of the following important Theorem \ref{hopf2}.
\begin{theorem}\label{hopf2}
Assume that each \(T_{V_i}\) is a finite-dimensional minuscule irreducible representation of \(U_{q}(\mathfrak{g}_i)\), for \(i=1,\dots,n\). Then

$(1)$ \
$(R_{V_1\otimes\cdots\otimes V_n})^{t_2}$ is invertible,
and $(R_{V_1\otimes\cdots\otimes V_n}^{-1})^{t_1}P(R_{V_1\otimes\cdots\otimes V_n}^{t_2})^{-1}PK_0=\textrm{const}\cdot K_0$;

$(2)$ \
$H_{R_{V_1\otimes\cdots\otimes V_n}}$ is a coquasitriangular Hopf subalgebra of $H_{R_{V_1}}\otimes \cdots\otimes H_{R_{V_n}}$.
\end{theorem}
\begin{proof}
The generators of the algebra  $H_{R_{V_1\otimes\cdots\otimes V_n}}$ are $1$,
$t^{\underline{i}}_{\underline{j}}$,
$\tilde{t}^{\underline{i}}_{\underline{j}}$,
which satisfy
\begin{gather}
R\underline{T}_{1}\underline{T}_{2}=\underline{T}_{2}\underline{T}_{1}R,\quad
R^{t}\widetilde{\underline{T}}_{1}\widetilde{\underline{T}}_{2}=\widetilde{\underline{T}}_{2}\widetilde{\underline{T}}_{1}R^{t},\quad
(R^{t_{2}})^{-1}\underline{T}_{1}\widetilde{\underline{T}}_{2}=\widetilde{\underline{T}}_{2}\underline{T}_{1}(R^{t_{2}})^{-1}.\label{tilde}
\end{gather}
Here,
$\underline{T}=(t^{\underline{i}}_{\underline{j}})$,
$\tilde{\underline{T}}=(t^{\underline{i}}_{\underline{j}})$
are $(m_1\cdots m_n)\times(m_1\cdots m_n)$-matrices.
The first assertion guarantees that the algebra $H_{R_{V_1\otimes\cdots\otimes V_n}}$ is a Hopf algebra.
We know that the Hopf algebra $H_{R_{V_1}}\otimes \cdots\otimes H_{R_{V_n}}$ is coquasitriangular, as each \(H_{R_{V_i}}\) is coquasitriangular.
There exists a natural map
\begin{equation}
\iota: H_{R_{V_1\otimes\cdots\otimes V_n}} \longrightarrow H_{R_{V_1}}\otimes\cdots\otimes H_{R_{V_n}}
\end{equation}
such that $t^{\underline{i}}_{\underline{j}} \rightarrow{t_{V_1}}^{i_1}_{j_1}\otimes\cdots\otimes {t_{V_n}}^{i_n}_{j_n}$
and $\tilde{t}^{\underline{i}}_{\underline{j}} \rightarrow \tilde{t}_{V_1}{}^{i_1}_{j_1}\otimes \cdots\otimes\tilde{t}_{V_n}{}^{i_n}_{j_n}$
if $\underline{i}=(i_{1},i_2,\cdots,i_n)$ and
$\underline{j}=(j_{1},j_2,\cdots,j_n)$,
where
${t_{V_k}}^{i_k}_{j_k}$ are the generators of $H_{R_{V_k}}$,
$k=1, \cdots, n$.

We will prove $\iota$ is a Hopf algebra embedding.

In terms of the algebra structure \eqref{tilde} of each $H_{R_{V_m}}$,
$m=1, \cdots, n$,
we have
\begin{gather}
\sum\limits_{a_m,b_m}(R_{V_m})^{j_mi_m}_{a_mb_m}(t_{V_m})^{a_m}_{k_m}(t_{V_m})^{b_m}_{l_m}
=(t_{V_m})^{i_m}_{a_m}(t_{V_m})^{j_m}_{b_m}(R_{V_m})^{b_ma_m}_{k_ml_m},\label{em1}\\
\sum\limits_{a_m,b_m}(R_{V_m}^{t})^{j_mi_m}_{a_mb_m}(\tilde{t}_{V_m})^{a_m}_{k_m}(\tilde{t}_{V_m})^{b_m}_{l_m}
=(\tilde{t}_{V_m})^{i_m}_{a_m}(\tilde{t}_{V_m})^{j_m}_{b_m}(R_{V_m}^{t})^{b_ma_m}_{k_ml_m},\label{em2}\\
\sum\limits_{a_m,b_m}((R_{V_m}^{t_2})^{-1})^{j_mi_m}_{a_mb_m}(t_{V_m})^{a_m}_{k_m}(\tilde{t}_{V_m})^{b_m}_{l_m}
=(\tilde{t}_{V_m})^{i_m}_{a_m}(t_{V_m})^{j_m}_{b_m}((R_{V_m}^{t_2})^{-1})^{b_ma_m}_{k_ml_m}\label{em3}.
\end{gather}
Then, using the algebra structure of \(H_{R_{V_1}}\otimes \cdots\otimes H_{R_{V_n}}\),
we derive the following relation inside \(H_{R_{V_1}}\otimes \cdots\otimes H_{R_{V_n}}\) via equality \eqref{em1}:
\begin{equation}\label{em4}
\left.
\begin{array}{rl}
&\sum\limits_{\underline{a},\underline{b}}(R_{V_1\otimes\cdots\otimes V_n})^{\underline{j}\,\underline{i}}_{\underline{a}\,\underline{b}}
\Bigl((t_{V_1})^{a_1}_{k_1}\otimes\cdots\otimes (t_{V_n})^{a_n}_{k_n}\Bigr)\Bigl((t_{V_1})^{b_1}_{l_1}\otimes\cdots\otimes(t_{V_n})^{b_n}_{l_n}\Bigr)\\
&\quad =\sum\limits_{\underline{a},\underline{b}}(R_{V_1})^{j_1i_1}_{a_1b_1}\cdots(R_{V_m})^{j_mi_m}_{a_mb_m}\Bigl((t_{V_1})^{a_1}_{k_1}\otimes\cdots\otimes
(t_{V_n})^{a_n}_{k_n}\Bigr)\Bigl((t_{V_1})^{b_1}_{l_1}\otimes\cdots\otimes(t_{V_n})^{b_n}_{l_n}\Bigr)\\
&\quad =\sum\limits_{\underline{a},\underline{b}}(R_{V_1})^{j_1i_1}_{a_1b_1}\cdots(R_{V_m})^{j_mi_m}_{a_mb_m}
\Bigl((t_{V_1})^{a_1}_{k_1}(t_{V_1})^{b_1}_{l_1}\otimes\cdots\otimes (t_{V_n})^{a_n}_{k_n}(t_{V_n})^{b_n}_{l_n}\Bigr)\\
&\quad =\sum\limits_{\underline{a},\underline{b}}\Bigl((t_{V_1})^{i_1}_{a_1}(t_{V_1})^{j_1}_{b_1}\otimes\cdots\otimes
(t_{V_n})^{i_n}_{a_n}(t_{V_n})^{j_n}_{b_n}\Bigr)(R_{V_1})^{b_1a_1}_{k_1l_1}\cdots(R_{V_n})^{b_na_n}_{k_nl_n},\\
&\quad =\sum\limits_{\underline{a},\underline{b}}\Bigl((t_{V_1})^{i_1}_{a_1}\otimes\cdots\otimes(t_{V_n})^{i_n}_{a_n}\Bigr)
\Bigl((t_{V_1})^{j_1}_{b_1}\otimes\cdots\otimes(t_{V_n})^{j_n}_{b_n}\Bigr)(R_{V_1})^{b_1a_1}_{k_1l_1}\cdots(R_{V_n})^{b_na_n}_{k_nl_n}\\
&\quad =\sum\limits_{\underline{a},\underline{b}}\Bigl((t_{V_1})^{i_1}_{a_1}\otimes\cdots\otimes(t_{V_n})^{i_n}_{a_n}\Bigr)
\Bigl((t_{V_1})^{j_1}_{b_1}\otimes\cdots\otimes(t_{V_n})^{j_n}_{b_n}\Bigr)(R_{V_1\otimes\cdots\otimes V_n})^{\underline{b}\,\underline{a}}_{\underline{k}\,\underline{l}}.
\end{array}
\right.
\end{equation}
From \eqref{em2} and \eqref{em3}, we obtain
\begin{equation}\label{em5}
\left.
\begin{array}{rl}
&\sum\limits_{\underline{a},\underline{b}}(R_{V_1\otimes\cdots\otimes V_n}^{t})^{\underline{j}\,\underline{i}}_{\underline{a}\,\underline{b}}
\Bigl((\tilde{t}_{V_1})^{a_1}_{k_1}\otimes\cdots\otimes (\tilde{t}_{V_n})^{a_n}_{k_n}\Bigr)\Bigl((\tilde{t}_{V_1})^{b_1}_{l_1}\otimes\cdots\otimes(\tilde{t}_{V_n})^{b_n}_{l_n}\Bigr)\\
&\quad =\sum\limits_{\underline{a},\underline{b}}\Bigl((\tilde{t}_{V_1})^{i_1}_{a_1}\otimes\cdots\otimes(\tilde{t}_{V_n})^{i_n}_{a_n}\Bigr)
\Bigl((\tilde{t}_{V_1})^{j_1}_{b_1}\otimes\cdots\otimes(\tilde{t}_{V_n})^{j_n}_{b_n}\Bigr)(R_{V_1\otimes\cdots\otimes
V_n}^{t})^{\underline{b}\,\underline{a}}_{\underline{k}\,\underline{l}},
\end{array}
\right.
\end{equation}
\begin{equation}\label{em6}
\left.
\begin{array}{rl}
&\sum\limits_{\underline{a},\underline{b}}((R_{V_1\otimes\cdots\otimes V_n}^{t_2})^{-1})^{\underline{j}\,\underline{i}}_{\underline{a}\,\underline{b}}
\Bigl((t_{V_1})^{a_1}_{k_1}\otimes\cdots\otimes (t_{V_n})^{a_n}_{k_n}\Bigr)\Bigl((\tilde{t}_{V_1})^{b_1}_{l_1}\otimes\cdots\otimes(\tilde{t}_{V_n})^{b_n}_{l_n}\Bigr)\\
&\quad =\sum\limits_{\underline{a},\underline{b}}\Bigl((t_{V_1})^{i_1}_{a_1}\otimes\cdots\otimes(t_{V_n})^{i_n}_{a_n}\Bigr)
\Bigl((\tilde{t}_{V_1})^{j_1}_{b_1}\otimes\cdots\otimes(\tilde{t}_{V_n})^{j_n}_{b_n}\Bigr)((R_{V_1\otimes\cdots\otimes
V_n}^{t_2})^{-1})^{\underline{b}\,\underline{a}}_{\underline{k}\,\underline{l}}.
\end{array}
\right.
\end{equation}

Relations \eqref{em4} --- \eqref{em6} imply that \(\iota\) is an algebra morphism.
Furthermore, it is easy to verify that
\begin{equation}
\Delta_{H_{R_{V_1}}\otimes \cdots\otimes H_{R_{V_n}}}\circ\iota=(\iota\otimes\iota)\circ\Delta_{H_{R_{V_1\otimes\cdots\otimes V_n}}},
\quad
\iota\circ S_{H_{R_{V_1\otimes\cdots\otimes V_n}}}=S_{H_{R_{V_1}}\otimes \cdots\otimes H_{R_{V_n}}}\circ\iota.
\end{equation}
Thus $\iota$ is a Hopf algebra homomorphism.

On the other hand, the Hopf algebra $H_{R_{V_1}}\otimes\cdots\otimes H_{R_{V_n}}$ is coquasitriangular. Its coquasitriangular structure is determined by those of the tensor-factor Hopf algebras \(H_{R_{V_i}}\), as described in Lemma 3.2, and the displayed relation yields the induced coquasitriangular structure on the subalgebra under consideration.
Accordingly, from the identity
\begin{gather*}
	\mathfrak{r}_{H_{R_{V_1\otimes\cdots\otimes V_n}}}(t^{\underline{i}}_{\underline{j}}\otimes t^{\underline{k}}_{\underline{l}})
	=(R_{V_1\otimes\cdots\otimes V_n})^{\underline{i}\,\underline{k}}_{\underline{j}\,\underline{l}},\qquad
	\mathfrak{r}_{H_{R_{V_1\otimes\cdots\otimes V_n}}}(t^{\underline{i}}_{\underline{j}}\otimes \tilde{t}^{\underline{k}}_{\underline{l}})
	=((R_{V_1\otimes\cdots\otimes V_n}^{t_2})^{-1})^{\underline{i}\,\underline{k}}_{\underline{j}\,\underline{l}},\\
	\mathfrak{r}_{H_{R_{V_1\otimes\cdots\otimes V_n}}}(\tilde{t}^{\underline{i}}_{\underline{j}}\otimes t^{\underline{k}}_{\underline{l}})
	=(R_{V_1\otimes\cdots\otimes V_n})^{-1}{}^{\underline{j}\,\underline{k}}_{\underline{i}\,\underline{l}},\qquad
	\mathfrak{r}_{H_{R_{V_1\otimes\cdots\otimes V_n}}}(\tilde{t}^{\underline{i}}_{\underline{j}}\otimes \tilde{t}^{\underline{k}}_{\underline{l}})
	=(R_{V_1\otimes\cdots\otimes V_n})^{\underline{j}\,\underline{l}}_{\underline{i}\,\underline{k}}.
\end{gather*}

This completes the proof.
\end{proof}
By Theorem \ref{hopf2} (2) and the pair $\left(U_q^{\text{ext}}(\mathfrak{g}_1)\otimes\cdots\otimes U_q^{\text{ext}}(\mathfrak{g}_n),H_{R_{V_1\otimes\cdots\otimes V_n}}\right)$ from Lemma \ref{weak},
we obtain
\begin{theorem}\label{weakly22}
There exists a weakly quasitriangular dual pair between
the quasitriangular Hopf algebra $U_q^{\text{ext}}(\mathfrak{g}_1\oplus\cdots\oplus\mathfrak{g}_n)$ and the coquasitriangular Hopf algebra
$H_{R_{V_1\otimes\cdots\otimes V_n}}$,
that is,
\begin{gather}
\langle (m^{+})^{\underline{i}}_{\underline{j}},t^{\underline{k}}_{\underline{l}}\rangle
=(R_{V_1\otimes\cdots\otimes V_n})^{\underline{i}\,\underline{k}}_{\underline{j},\underline{l}},\qquad
\langle (m^{-})^{\underline{i}}_{\underline{j}},t^{\underline{k}}_{\underline{l}}\rangle=(R^{-1}_{V_1\otimes\cdots\otimes
V_n}){}^{\underline{k}\,\underline{i}}_{\underline{l}\,\underline{j}},\label{tlidew1}\\
\langle (m^{+})^{\underline{i}}_{\underline{j}},\tilde{t}^{\underline{k}}_{\underline{l}}\rangle=(R_{V_1\otimes\cdots\otimes
V_n}^{t_{2}})^{-1}{}^{\underline{i}\,\underline{k}}_{\underline{j}\,\underline{l}},\qquad
\langle (m^{-})^{\underline{i}}_{\underline{j}},\tilde{t}^{\underline{k}}_{\underline{l}}\rangle=[(R^{-1}_{V_1\otimes\cdots\otimes
V_n})^{t_{1}}]^{-1}{}^{\underline{k}\,\underline{i}}_{\underline{l}\,\underline{j}}.\label{tlidew2}
\end{gather}
The convolution-invertible algebra\,/\,anti-coalgebra maps $\mathcal{R},\bar{\mathcal{R}}$ are
$$
\mathcal{R}(t^{\underline{i}}_{\underline{j}})=(m^{+})^{\underline{i}}_{\underline{j}},\qquad
\mathcal{R}(\tilde{t}^{\underline{i}}_{\underline{j}})=(m^{+})^{-1}{}^{\underline{j}}_{\underline{i}};\qquad
\bar{\mathcal{R}}(t^{\underline{i}}_{\underline{j}})=(m^{-})^{\underline{i}}_{\underline{j}},\qquad
\bar{\mathcal{R}}(\tilde{t}^{\underline{i}}_{\underline{j}})=(m^{-})^{-1}{}^{\underline{j}}_{\underline{i}}.
$$
\end{theorem}
Theorem \ref{weakly22} may also be proved directly by an argument analogous to that given in \cite{HH2}.

\subsection{Multi-tensor-product of generalized double-bosonization}
Using the weakly quasi-triangular dual pair from Theorem \ref{weakly22}, one strategy for constructing dually-paired braided groups within the braided module categories \(\mathfrak{M}_{U_q^{\text{ext}}(\mathfrak{g}_1\oplus\cdots\oplus\mathfrak{g}_n)}\) and \(_{U_q^{\text{ext}}(\mathfrak{g}_1\oplus\cdots\oplus\mathfrak{g}_n)}\mathfrak{M}\) is to locate the corresponding objects in the braided comodule categories \(^{H_{R_{V_1\otimes\cdots\otimes V_n}}}\mathfrak{M}\) and \(\mathfrak{M}^{H_{R_{V_1\otimes\cdots\otimes V_n}}}\).
The pair \(V(R',R)\) and \(V^{\vee}(R',R_{21}^{-1})\) are then dually-paired braided groups belonging to the braided category of right- and left- \(U_q^{\text{ext}}(\mathfrak{g}_1\oplus\cdots\oplus\mathfrak{g}_n)\)-modules, respectively.

To maintain consistent notation, denote by $e^{\underline{i}}, f_{\underline{i}}$ the generators of $V(R',R), V^{\vee}(R',R_{21}^{-1})$ associated with the $R$-matrix \(R_{V_1\otimes\cdots\otimes V_n}\), respectively,
where \(\underline{i}=(i_1,\cdots,i_n)\) with \(i_k=1,\cdots,\) \(m_k\) for \(k=1,\cdots,n\). 
In the double-bosonization construction, one needs an additional new group-like element in the Cartan sector, which arises from adjoining the connecting node as the common grafting vertex for several sub-Dynkin diagrams.
Accordingly, we must refine the weakly quasitriangular dual pair from Theorem \ref{weakly22} to its central extensions
$$
\Bigl(\widetilde{U_{q}^{\textrm{ext}}}(\mathfrak{g}_1\oplus\cdots\oplus\mathfrak{g}_n),\; \widetilde{H}_{R_{V_1\otimes\cdots\otimes V_n}}\Bigr),
$$
where
\(\widetilde{U_{q}^{\textrm{ext}}}(\mathfrak{g}_1\oplus\cdots\oplus\mathfrak{g}_n)=U_{q}^{\textrm{ext}}(\mathfrak{g}_1\oplus\cdots\oplus\mathfrak{g}_n)\otimes k[c,c^{-1}]\),
\(\widetilde{H}_{R_{V_1\otimes\cdots\otimes V_n}}=H_{R_{V_1\otimes\cdots\otimes V_n}}\otimes k[g,g^{-1}]\),
subject to \(\langle c,g\rangle=\lambda\), \(\mathscr{R}(g)=c^{-1}\), and \(\bar{\mathscr{R}}(g)=c\), with \(R_{V_1\otimes\cdots\otimes V_n}=\lambda R\).
Then our concerned braided objects $V(R^{\prime},R)\in
{}^{\widetilde{H}_{R_{V_1\otimes\cdots\otimes V_n}}}\mathfrak M$,
$V^{\vee}(R^{\prime},R_{21}^{-1})\in
 \mathfrak M^{\widetilde{H}_{R_{V_1\otimes\cdots\otimes V_n}}}$ with the
coactions $e^{\underline{i}} \mapsto g\,t^{\underline{i}}_{\underline{a}}\otimes e^{\underline{a}}$ and
$ f_{\underline{i}} \mapsto f_{\underline{a}}\otimes g\,t^{\underline{i}}_{\underline{a}}$.
Thus we have dually-paired braided groups $V(R',R)$ and $V^{\vee}(R',R_{21}^{-1})$ in the braided  $\widetilde{U_q^{\text{ext}}}(\mathfrak{g}_1\oplus\cdots\oplus\mathfrak{g}_n)$-modules category.

With these points in mind, we arrive at the following {\it grafting-construction theorem}, which can be viewed as the multi-tensor-product version of generalized double-bosonization construction stated in Theorem \ref{cor1}.
Let \(T_{V_i}\) denote a minuscule irreducible representation of \(U_q^{\text{ext}}(\mathfrak{g}_i)\) with minuscule highest weight (see \cite{FH}).
\begin{theorem}\label{multi}
Let $R_{V_1\otimes\cdots\otimes V_n}$ be the $R$-matrix of the representation $T_{V_1\otimes\cdots\otimes V_n}$ of $U_q^{\text{ext}}(\mathfrak{g}_1\oplus\cdots\oplus\mathfrak{g}_n)$, and $\lambda$ a non-zero scalar such that the Majid pair $(R=\lambda^{-1}R_{V_1\otimes\cdots\otimes V_n}, R')$ satisfies conditions \eqref{*} and \eqref{**}.
Then the tensor vector space $V^{\vee}(R^{\prime},R_{21}^{-1})\otimes \widetilde{U_{q}^{\textrm{ext}}}(\mathfrak{g}_1\oplus\cdots\oplus\mathfrak{g}_n)\otimes V(R^{\prime},R)$ carries a Hopf-algebra structure inherited from double-bosonization, denoted
$$
U=U\left(V^{\vee}(R^{\prime},R_{21}^{-1}),\widetilde{U_{q}^{\textrm{ext}}}(\mathfrak{g}_1\oplus\cdots\oplus\mathfrak{g}_n),V(R^{\prime},R)\right).
$$
Its algebra structure is governed by the cross relations
\begin{gather*}
cf_{\underline{i}}=\lambda f_{\underline{i}}c,\quad
e^{\underline{i}}c=\lambda ce^{\underline{i}},\quad
[c,(m^{\pm})^{\underline{i}}_{\underline{j}}]=0,\quad
[e^{\underline{i}},f_{\underline{j}}]=
\frac{(m^{+})^{\underline{i}}_{\underline{j}}c^{-1}-c(m^{-})^{\underline{i}}_{\underline{}j}}{q_{\ast}-q_{\ast}^{-1}},\\
e^{\underline{i}}(m^{+})^{\underline{j}}_{\underline{k}}
=R_{VV}{}^{\underline{j}\,\underline{i}}_{\underline{a}\,\underline{b}}(m^{+})^{\underline{a}}_{\underline{k}}e^{\underline{b}},\quad
(m^{-})^{\underline{i}}_{\underline{j}}e^{\underline{k}}
=R_{VV}{}^{\underline{k}\,\underline{i}}_{\underline{a}\,\underline{b}}e^{\underline{a}}(m^{-})^{\underline{b}}_{\underline{j}},\\
(m^{+})^{\underline{i}}_{\underline{j}}f_{\underline{k}}=f_{\underline{b}}(m^{+})^{\underline{i}}_{\underline{a}}
R_{VV}{}^{\underline{a}\,\underline{b}}_{\underline{j}\,\underline{k}},\quad
f_{\underline{i}}(m^{-})^{\underline{j}}_{\underline{k}}=(m^{-})^{\underline{j}}_{\underline{b}}f_{\underline{a}}
R_{VV}{}^{\underline{a}\,\underline{b}}_{\underline{i}\,\underline{k}},
\end{gather*}
while its coalgebra structure is given by the comultiplication and counit
\begin{gather*}
\Delta (c)=c\otimes c, \quad \Delta
(e^{\underline{i}})=e^{\underline{a}}\otimes (m^{+})^{\underline{i}}_{\underline{a}}c^{-1}+1\otimes e^{\underline{i}}, \quad
\Delta (f_{\underline{i}})=f_{\underline{i}}\otimes 1+c(m^{-})^{\underline{a}}_{\underline{i}}\otimes f_{\underline{a}},
\\
\epsilon (e^{\underline{i}})=\epsilon (f_{\underline{i}})=0.
\end{gather*}
Here one may normalize \(e^{\underline{i}}\) so that \(q_{\ast}-q_{\ast}^{-1}\) satisfies the standard convention.
\end{theorem}

\section{Applications: Grafting Constructions of Two Quantum Groups $U_q(\mathfrak{sl}_{m+n})$ and $U_q(F_4)$}
We focus on the case \(n=2\) of Theorem \ref{multi}, with the aim of describing the grafting construction for certain quantum enveloping algebras. As usual, let
\(\mathfrak{sl}_n\) denote the finite-dimensional complex simple Lie algebra of type \(A_{n-1}\), with simple roots
\(\alpha_i\) satisfying \((\alpha_i,\alpha_i)=2\).
Let \(\epsilon_i\) stand for the fundamental weight corresponding to
the simple root \(\alpha_i\), for \(i=1,\dots,n-1\).
As is well known,
the natural representation $T_{\mathbb{C}^{n}}$ of $U_q(\mathfrak{sl}_{n})$ is
\begin{equation}\label{vector}
T_{\mathbb{C}^{n}}(E_i)(v_{i+1})=v_i, \
T_{\mathbb{C}^{n}}(F_i)(v_i)=v_{i+1}, \
T_{\mathbb{C}^{n}}(K_i)(v_i)=qv_{i}, \
T_{\mathbb{C}^{n}}(K_i)(v_{i+1})=q^{-1}v_{i+1},
\end{equation}
where $v_1,\cdots,v_n$ is a basis of the $n$-dimensional vector space $\mathbb{C}^{n}$.
The corresponding $R$-matrix $R_{\mathbb{C}^{n}}$ has the form
\begin{equation}\label{matrix}
R_{\mathbb{C}^{n}}=q\sum\limits_{i=1}^{n}E_{ii}\otimes E_{ii}+\sum_{i,j=1,i\neq j}^{n}E_{ii}\otimes E_{jj}+(q-q^{-1})
\sum_{i,j=1,i> j}^{n}E_{ij}\otimes E_{ji},
\end{equation}
where $E_{ij}$ is the $n\times n$ matrix with $1$ in the $(i,j)$-position and $0$ elsewhere.

The representation dual to $T_{\mathbb{C}^{n}}$ for $U_q(\mathfrak{sl}_{n})$ is
$T_{(\mathbb{C}^{n})^*}(E_i)(v_{i}^{*})=v_{i+1}^{*}$,
$T_{(\mathbb{C}^{n})^*}(F_i)(v_{i+1}^{*})=v_{i}^{*}$,
$T_{(\mathbb{C}^{n})^*}(K_i)(v_i^{*})=q^{-1}v_{i}^{*}$,
$T_{(\mathbb{C}^{n})^*}(K_i)(v_{i+1}^{*})=qv_{i+1}^{*}$.
The upper/lower triangular form of $m^+$/$m^-$ follows from the lowering/raising action of $T_{(\mathbb C^n)^*}(E_i)$ and $T_{(\mathbb C^n)^*}(F_i)$, respectively.
For a linear transformation, its matrices in two oppositely ordered bases are mutual transposes. Hence we reverse the indices of the basis of the original dual representation via setting $v_i^{*}$ by $\omega_{n+1-i}$,
and then we have that
\begin{equation}\label{dual}
\left.
\begin{array}{c}
T_{(\mathbb{C}^{n})^*}(E_{n-i})(\omega_{i{+}1})=\omega_{i},\quad
T_{(\mathbb{C}^{n})^*}(F_{n-i})(\omega_{i})=\omega_{i+1},\\
T_{(\mathbb{C}^{n})^*}(K_{n-i})(\omega_{i})=q\omega_{i},\quad
T_{(\mathbb{C}^{n})^*}(K_{n-i})(\omega_{i{+}1})=q^{-1}\omega_{i+1}.
\end{array}
\right.
\end{equation}
By using equations \eqref{vector}, \eqref{dual}, together with the stated weight calculation, the $R$-matrix corresponding to $T_{V^{*}}$ is the one defined in \eqref{matrix}.
Then each entry in the $R$-matrix is $$(R_{\mathbb{C}^{n}})^{ij}_{kl}=q^{\delta_{ij}}\delta_{ik}\delta_{jl}+(q-q^{-1})\delta_{il}\delta_{jk}\theta(j-i),$$
and the braiding $PR_{\mathbb{C}^{n}}$ satisfies $$(PR_{\mathbb{C}^{n}}-q\textrm{I})(PR_{\mathbb{C}^{n}}+q^{-1}\textrm{I})=0,
\mbox{\quad where\quad} \theta(k)= \left\{
\begin{array}{lcl}
1&~~&k>0,\\
0&~~&k\leq 0.
\end{array}
\right.
$$

\subsection{Simply-laced case: Type $A$}
In this subsection,
we aim to obtain $U_q(\mathfrak{sl}_{n+m})$ by grafting on the vector space $V^{\vee}(R^{\prime},R_{21}^{-1})\otimes (U_q^{\text{ext}}(\mathfrak{sl}_n\oplus\mathfrak{sl}_m)\otimes
k[c,c^{-1}])\otimes V(R^{\prime},R)$.
To do this,
we choose the vector representation $\mathbb{C}^n$ of $U_q(\mathfrak{sl}_{n})$
and  the dual vector representation $(\mathbb{C}^m)^*$ of $U_q(\mathfrak{sl}_{m})$.
By Proposition~\ref{minimal}, the braiding
$PR_{\mathbb{C}^n\otimes(\mathbb{C}^m)^*}$ is diagonalizable with
distinct eigenvalues $q^2$, $-1$, and $q^{-2}$, hence its minimal
polynomial equation is
\begin{equation}\label{new}
(PR_{\mathbb{C}^n\otimes(\mathbb{C}^m)^*}+\textrm{I})(PR_{\mathbb{C}^n\otimes(\mathbb{C}^m)^*}-q^{2}\textrm{I})(PR_{\mathbb{C}^n\otimes(\mathbb{C}^m)^*}-q^{-2}\textrm{I})=0.
\end{equation}
We normalize $R_{\mathbb{C}^n\otimes(\mathbb{C}^m)^*}$ at the eigenvalue $-1$, that is, set
\begin{equation}\label{rmatrix1}
R=R_{\mathbb{C}^n\otimes(\mathbb{C}^m)^*},\quad
R'=RPR-(q^{2}+q^{-2})R+2P,
\end{equation}
so that $(PR+I)(PR'-I)=0$.

Furthermore, using the resulting dually paired braided groups \(V^{\vee}(R^{\prime},R_{21}^{-1})\) and \(V(R',R)\),
we arrive at the following grafting result.
\begin{theorem}\label{typeA}
Upon identifying the braided-group generators $e^{(n,m)}$ and
$f_{(n,m)}$, together with the group-like element
$(m^+)^{(n,m)}_{(n,m)}c^{-1}$, with the extra simple-root generators \(E_n\), \(F_n\), and \(K_n\) respectively,
the resulting quantum group
$$
U=U\Bigl(V^{\vee}(R^{\prime},R_{21}^{-1}),U_q^{\text{ext}}(\mathfrak{sl}_n\oplus\mathfrak{sl}_m)\otimes k[c,c^{-1}],V(R^{\prime},R)\Bigr)
$$
is isomorphic to the quantum enveloping algebra $U_q(\mathfrak{sl}_{n+m})$ with $K_i^{\pm\frac{1}{n}}$,
$i=1, \cdots, n-1$
and
$K_j^{\pm\frac{1}{m}}$,
$j=n+1, \cdots, n+m-1$ adjoined.
The extra simple root is represented by the filled vertex in the Dynkin diagram below:
\begin{center}
\setlength{\unitlength}{1mm}
\begin{picture}(68,8)
\put(2,4){\circle{1}}
\put(0,6){$\alpha_1$}
\put(2.5,4){\line(1,0){6}}
\put(9,4){\circle{1}}
\put(7,6){$\alpha_2$}
\put(9.5,4){\line(1,0){6}}
\put(16,4){$\ldots$}
\put(20.5,4){\line(1,0){6}}
\put(27,4){\circle{1}}
\put(24,6){$\alpha_{n-1}$}
\multiput(27.5,4)(1,0){6}{\line(1,0){0.5}}
\put(34,4){\circle*{1}}
\multiput(34.5,4)(1,0){6}{\line(1,0){0.5}}
\put(10,0){$\mathfrak{g}_1\simeq \mathfrak{sl}_{n}$}
\put(51,0){$\mathfrak{g}_2\simeq \mathfrak{sl}_{m}$}
\put(41,4){\circle{1}}
\put(38,6){$\alpha_{n+1}$}
\put(41.5,4){\line(1,0){6}}
\put(48,4){\circle{1}}
\put(46,6){$\alpha_{n+2}$}
\put(48.5,4){\line(1,0){6}}
\put(55,4){$\ldots$}
\put(59.5,4){\line(1,0){6}}
\put(66,4){\circle{1}}
\put(63,6){$\alpha_{n+m-1}$}
\end{picture}
\end{center}
\end{theorem}
\begin{proof}
Let $\{E_i, F_i, K_i\mid i=1,\cdots, n-1\}$ denote the generators of the quantum enveloping algebra $U_q(\mathfrak{sl}_n)$, and let $\{E_i, F_i, K_i\mid i=n+1,\cdots, n+m-1\}$ denote the generators of $U_q(\mathfrak{sl}_m)$, the FRT-generators $m^{\pm}$ associated with the vector representation $\mathbb{C}^n$ are given below (see \cite{HH1}):
\begin{equation}\label{A0}
\left.
\begin{array}{c}
(m_{\mathbb{C}^n}^{+})^{i}_{i}=K^{-\frac{1}{n}}_{1}K^{-\frac{2}{n}}_{2}\cdots K^{-\frac{i-1}{n}}_{i-1}
K^{\frac{n-i}{n}}_{i}\cdots K^{\frac{n-(n-1)}{n}}_{n-1},\qquad
\quad (m_{\mathbb{C}^n}^{+})^{i}_{i}=K_i(m_{\mathbb{C}^n}^{+})^{i+1}_{i+1},\\
(m_{\mathbb{C}^n}^{+})^{i}_{i+1}=(q-q^{-1})E_{i}(m_{\mathbb{C}^n}^{+})^{i+1}_{i+1}, \qquad\quad
(m_{\mathbb{C}^n}^{-})^{i+1}_{i}=(q-q^{-1})(m_{\mathbb{C}^n}^{-})^{i+1}_{i+1}F_{i}.
\end{array}
\right.
\end{equation}
For the dual vector representation $T_{(\mathbb{C}^m)^*}$ of $U_q(\mathfrak{sl}_m)$,
a direct calculation gives
\begin{equation}\label{A0'}
\left.
\begin{array}{c}
(m_{(\mathbb{C}^m)^*}^{+})^{i}_{i}=K^{-\frac{1}{m}}_{n+m-1}{\cdots}
K^{\frac{m-i}{m}}_{n+m-i}{\cdots} K^{\frac{m-(m-1)}{m}}_{n+m-(m-1)},\quad
(m_{(\mathbb{C}^m)^*}^{+})^{i}_{i}=K_{n+m-i}(m_{(\mathbb{C}^m)^*}^{+})^{i+1}_{i+1},\\
(m_{(\mathbb{C}^m)^*}^{+})^{i}_{i+1}=(q{-}q^{-1})E_{n+m-i}(m_{(\mathbb{C}^m)^*}^{+})^{i+1}_{i+1}, \quad
(m_{(\mathbb{C}^m)^*}^{-})^{i+1}_{i}=(q{-}q^{-1})(m_{(\mathbb{C}^m)^*}^{-})^{i+1}_{i+1}F_{n+m-i}.
\end{array}
\right.
\end{equation}
These FRT-generators $m^{\pm}$ generate the classical extended Hopf algebra $U_q^{\textrm{ext}}(\mathfrak{sl}_n\oplus \mathfrak{sl}_m)$.

\emph{\textbf{Step 1}. The cross relations between the new simple root vector $E_n$ and $K_i$'s.}

Under this identification, and using the cross relations from Theorem \ref{multi}, we obtain
\begin{equation}\label{A1}
\left.
\begin{array}{rl}
E_nK_n&=e^{(n,m)}(m^+)^{(n,m)}_{(n,m)}c^{-1}=\lambda R^{(nm)(nm)}_{(ab)(cd)}(m^+)^{(a,b)}_{(n,m)}e^{(c,d)}c^{-1}\\
&=R^{(nm)(nm)}_{(ab)(cd)}(m^+)^{(a,b)}_{(n,m)}c^{-1}e^{(c,d)}=(R_{\mathbb{C}^n\otimes(\mathbb{C}^m)^*})^{(nm)(nm)}_{(ab)(cd)}(m^+)^{(a,b)}_{(n,m)}c^{-1}e^{(c,d)}\\
&=(R_{\mathbb{C}^n})^{nn}_{ac}(R_{(\mathbb{C}^m)^*})^{mm}_{bd}\Bigl((m_{\mathbb{C}^n}^{+})^a_n\otimes (m_{(\mathbb{C}^m)^*}^{+})^b_m\Bigr)\,c^{-1}e^{(c,d)}\\
&=(R_{\mathbb{C}^n})^{nn}_{nn}(R_{(\mathbb{C}^m)^*})^{mm}_{mm}\Bigl((m_{\mathbb{C}^n}^{+})^n_n\otimes (m_{(\mathbb{C}^m)^*}^{+})^m_m\Bigr)\,c^{-1}e^{(n,m)}\\
&=qq(m^+)^{(n,m)}_{(n,m)}\,c^{-1}e^{(n,m)}=q^{2}K_nE_n.
\end{array}
\right.
\end{equation}

We are in a position to determine the relations between $e^{(n,m)}$ and group-likes $K_i$.
The relations $(m_{\mathbb{C}^n}^{+})^{i}_{i}=K_i(m_{\mathbb{C}^n}^{+})^{i+1}_{i+1},\,1\leq i\leq n-1,$ and $(m_{(\mathbb{C}^m)^*}^{+})^{i}_{i}=K_{n+m-i}(m_{(\mathbb{C}^m)^*}^{+})^{i+1}_{i+1},\,
1\leq i\leq m{-}1$
imply the following equalities in quantum group
$U_q^{\text{ext}}(\mathfrak{sl}_n\oplus\mathfrak{sl}_m)$:
\begin{gather}
(m^+)^{(i,m)}_{(i,m)}
=(m_{\mathbb{C}^n}^{+})^{i}_{i}\otimes (m_{(\mathbb{C}^m)^*}^{+})^{m}_{m}
=(K_{i}\otimes 1)((m_{\mathbb{C}^n}^{+})^{i+1}_{i+1}\otimes(m_{(\mathbb{C}^m)^*}^{+})^{m}_{m})
=K_{i}(m^+)^{(i+1,m)}_{(i+1,m)}, \label{A2}\\
(m^+)^{(n,i)}_{(n,i)}
=(m_{\mathbb{C}^n}^{+})^{n}_{n}{\otimes} (m_{(\mathbb{C}^m)^*}^{+})^{i}_{i}
=(1{\otimes} K_{n+m-i})(m_{\mathbb{C}^n}^{+})^{n}_{n}{\otimes}(m_{(\mathbb{C}^m)^*}^{+})^{i+1}_{i+1}
=K_{n+m-i}(m^+)^{(n,i+1)}_{(n,i+1)}.\label{A3}
\end{gather}

By virtue of the cross relation $e^{\underline{i}}(m^{+})^{\underline{j}}_{\underline{k}}
=\lambda R^{\underline{j}\,\underline{i}}_{\underline{a}\,\underline{b}}(m^{+})^{\underline{a}}_{\underline{k}}e^{\underline{b}}$
from Theorem \ref{multi}, we obtain
\begin{equation}\label{A4}
\left.
\begin{array}{rl}
E_n\,(m^+)^{(j,m)}_{(j,m)}&=e^{(n,m)}\,(m^+)^{(j,m)}_{(j,m)}=\lambda R^{(jm)(nm)}_{(ab)(cd)}\,(m^+)^{(a,b)}_{(j,m)}\,e^{(c,d)}\\
&=(R_{\mathbb{C}^n})^{jn}_{ac}(R_{(\mathbb{C}^m)^*})^{mm}_{bd}\Bigl((m_{\mathbb{C}^n}^{+})^a_j\otimes (m_{(\mathbb{C}^m)^*}^{+})^b_m\Bigr)\,e^{(c,d)}\\
&=(R_{\mathbb{C}^n})^{jn}_{jn}(R_{(\mathbb{C}^m)^*})^{mm}_{mm}\Bigl((m_{\mathbb{C}^n}^{+})^j_j\otimes (m_{(\mathbb{C}^m)^*}^{+})^m_m\Bigr)\,e^{(n,m)}\\
&=(R_{\mathbb{C}^n})^{jn}_{jn}(R_{(\mathbb{C}^m)^*})^{mm}_{mm}(m^+)^{(j,m)}_{(j,m)}\,E_n,
\end{array}
\right.
\end{equation}
for any $1\leq j\leq n$.

Similarly, we have
\begin{equation}\label{A4'}
E_n\,(m^+)^{(n,k)}_{(n,k)}=(R_{\mathbb{C}^n})^{nn}_{nn}\,(R_{(\mathbb{C}^m)^*})^{km}_{km}\,(m^+)^{(n,k)}_{(n,k)}\,E_n,
\end{equation}
for all \(1\leq k\leq m\).

With equalities \eqref{A4} and \eqref{A4'} in hand, we obtain the following relations:
\begin{equation}\label{A5}
\left.
\begin{array}{rl}
E_nK_{i}(m^+)^{(i+1,m)}_{(i+1,m)}&=(R_{\mathbb{C}^n})^{in}_{in}(R_{(\mathbb{C}^m)^*})^{mm}_{mm}K_{i}(m^+)^{(i+1,m)}_{(i+1,m)}E_n\\
&\overset{\eqref{A4}}=(R_{\mathbb{C}^n})^{in}_{in}(R_{(\mathbb{C}^m)^*})^{mm}_{mm}K_{i}\frac{1}{(R_{\mathbb{C}^n})^{i+1,n}_{i+1,n}(R_{(\mathbb{C}^m)^*})^{mm}_{mm}}
E_n(m^+)^{(i+1,m)}_{(i+1,m)}\\
&=\frac{(R_{\mathbb{C}^n})^{in}_{in}}{(R_{\mathbb{C}^n})^{i+1,n}_{i+1,n}}K_{i}E_n(m^+)^{(i+1,m)}_{(i+1,m)},
\end{array}
\right.
\end{equation}
\begin{equation}\label{A66}
\left.
\begin{array}{rl}
E_nK_{n+m-i}(m^+)^{(n,i+1)}_{(n,i+1)}&=(R_{\mathbb{C}^n})^{nn}_{nn}(R_{(\mathbb{C}^m)^*})^{im}_{im}K_{n+m-i}(m^+)^{(n,i+1)}_{(n,i+1)}E_n\\
&\overset{\eqref{A4'}}=(R_{\mathbb{C}^n})^{nn}_{nn}(R_{(\mathbb{C}^m)^*})^{im}_{im}K_{n+m-i}\frac{1}{(R_{\mathbb{C}^n})^{nn}_{nn}(R_{(\mathbb{C}^m)^*})^{i+1,m}_{i+1,m}}
E_n(m^+)^{(n,i+1)}_{(n,i+1)}\\
&=\frac{(R_{(\mathbb{C}^m)^*})^{im}_{im}}{(R_{(\mathbb{C}^m)^*})^{i+1,m}_{i+1,m}}K_{n+m-i}E_n(m^+)^{(n,i+1)}_{(n,i+1)}
\end{array}
\right.
\end{equation}
in terms of \eqref{A2} and
\eqref{A3}.

Multiplying both sides of equalities \eqref{A5} and \eqref{A66} on the right by \((m^-)^{(i+1,m)}_{(i+1,m)}\) and \((m^-)^{(n,i+1)}_{(n,i+1)}\), respectively, we obtain
\begin{equation}\label{A7}
E_nK_{i}=\frac{(R_{\mathbb{C}^n})^{in}_{in}}{(R_{\mathbb{C}^n})^{i+1,n}_{i+1,n}}K_{i}E_n,\quad
E_nK_{n+m-i}=\frac{(R_{(\mathbb{C}^m)^*})^{im}_{im}}{(R_{(\mathbb{C}^m)^*})^{i+1,m}_{i+1,m}}K_{n+m-i}E_n.
\end{equation}
Then, corresponding to
$$R_{V}{}^{jk}_{jk}=q^{\delta_{jk}}+(q-q^{-1})\delta_{jk}\theta(k-j)=
\left\{
\begin{array}{ll}
1,&j\neq k;\\
q,&j=k
\end{array}
\right.
\mbox{\quad for any~} j, k,$$
we conclude that
\begin{equation}\label{A6}
\left.
\begin{array}{l}
E_nK_{i}=K_{i}E_n,\quad 1\leq i\leq n-2 \quad\mbox{and}\quad n+2\leq i\leq n+m-1,\\
E_nK_{n-1}=q^{-1}K_{n-1}E_n,\quad
E_nK_{n+1}=q^{-1}K_{n+1}E_n.
\end{array}
\right.
\end{equation}

\emph{\textbf{Step 2}. The cross relations between the new  additional group-like $K_n$ and $E_i$'s.}

Under the identification from Theorem \ref{typeA}, the extra group-like element \(K_n\) is given by $K^{{-}\frac{1}{n}}_{1}K^{{-}\frac{2}{n}}_{2}{\cdots}K^{{-}\frac{n{-}1}{n}}_{n{-}1}
K_{n+1}^{{-}\frac{m{-}1}{m}}K_{n{+}2}^{{-}\frac{m{-}2}{m}}{\cdots}K_{n{+}m{-}1}^{{-}\frac{1}{m}}c^{{-}1}$.
We then obtain the cross relations between \(K_n\) and \(E_i\) for \(1\leq i\leq n-1\) and \(n+1\leq i\leq n+m-1\), respectively:
\begin{equation*}
\begin{split}
E_{1}K_n&=E_1K^{-\frac{1}{n}}_{1}K^{-\frac{2}{n}}_{2}\cdots K^{-\frac{n-1}{n}}_{n-1}K_{n+1}^{-\frac{m-1}{m}}\cdots K_{n+m-1}^{-\frac{1}{m}}c^{-1}\\
&=q^{-\frac{2}{n}}q^{\frac{2}{n}}K^{-\frac{1}{n}}_{1}K^{-\frac{2}{n}}_{2}\cdots K^{-\frac{n-1}{n}}_{n-1}K_{n+1}^{-\frac{m-1}{m}}\cdots
K_{n+m-1}^{-\frac{1}{m}}c^{-1}E_1\\
&=K_nE_1,
\end{split}
\end{equation*}
\begin{equation*}
\begin{split}
E_{i}K_n&=E_iK^{-\frac{1}{n}}_{1}K^{-\frac{2}{n}}_{2}\cdots K^{-\frac{n-1}{n}}_{n-1}K_{n+1}^{-\frac{m-1}{m}}\cdots K_{n+m-1}^{-\frac{1}{m}}c^{-1}\\
&=q^{\frac{i-1}{n}}q^{-\frac{2i}{n}}q^{\frac{i+1}{n}}K^{-\frac{1}{n}}_{1}\cdots K^{-\frac{n-1}{n}}_{n-1}\cdots K_{n+1}^{-\frac{m-1}{m}}\cdots
K_{n+m-1}^{-\frac{1}{m}}c^{-1}E_i\\
&=K_nE_i,\quad 2\leq i\leq n-2,
\\
E_{n-1}K_n&=E_{n-1}K^{-\frac{1}{n}}_{1}K^{-\frac{2}{n}}_{2}\cdots K^{-\frac{n-1}{n}}_{n-1} K_{n+1}^{-\frac{m-1}{m}}\cdots K_{n+m-1}^{-\frac{1}{m}}c^{-1}\\
&=q^{\frac{n-2}{n}}q^{-\frac{2(n-1)}{n}}K^{-\frac{1}{n}}_{1}K^{-\frac{2}{n}}_{2}\cdots K^{-\frac{n-1}{n}}_{n-1} K_{n+1}^{-\frac{m-1}{m}}\cdots
K_{n+m-1}^{-\frac{1}{m}}c^{-1}E_{n-1}\\
&=q^{-1}K_nE_{n-1},
\\
E_{n+1}K_n&=E_{n+1}K^{-\frac{1}{n}}_{1}K^{-\frac{2}{n}}_{2}\cdots K^{-\frac{n-1}{n}}_{n-1} K_{n+1}^{-\frac{m-1}{m}}\cdots K_{n+m-1}^{-\frac{1}{m}}c^{-1}\\
&=K^{-\frac{1}{n}}_{1}K^{-\frac{2}{n}}_{2}\cdots K^{-\frac{n-1}{n}}_{n-1}E_{n+1}K_{n+1}^{-\frac{m-1}{m}}\cdots K_{n+m-1}^{-\frac{1}{m}}c^{-1}\\
&=q^{-\frac{2(m-1)}{m}}q^{\frac{m-2}{m}}K^{-\frac{1}{n}}_{1}K^{-\frac{2}{n}}_{2}\cdots K^{-\frac{n-1}{n}}_{n-1}K_{n+1}^{-\frac{m-1}{m}}\cdots
K_{n+m-1}^{-\frac{1}{m}}c^{-1}E_{n+1}\\
&=q^{-1}K_nE_{n+1},
\\
E_{j}K_n&=K^{-\frac{1}{n}}_{1}K^{-\frac{2}{n}}_{2}\cdots K^{-\frac{n-1}{n}}_{n-1}E_jK_{n+1}^{-\frac{m-1}{m}}\cdots K_{n+m-1}^{-\frac{1}{m}}c^{-1}\\
&=q^{\frac{n+m-(j-1)}{m}}q^{-\frac{2(n+m-j)}{m}}q^{\frac{n+m-(j+1)}{m}}K^{-\frac{1}{n}}_{1}\cdots K^{-\frac{n-1}{n}}_{n-1}\cdots K_{n+1}^{-\frac{m-1}{m}}\cdots
K_{n+m-1}^{-\frac{1}{m}}c^{-1}E_j\\
&=K_nE_j,\quad n+2\leq j\leq n+m-2,
\\
E_{n+m-1}K_n&=K^{-\frac{1}{n}}_{1}K^{-\frac{2}{n}}_{2}\cdots K^{-\frac{n-1}{n}}_{n-1} K_{n+1}^{-\frac{m-1}{m}}\cdots
E_{n+m-1}K_{n+m-2}^{-\frac{2}{m}}K_{n+m-1}^{-\frac{1}{m}}c^{-1}\\
&=q^{\frac{2}{m}}q^{-2\frac{1}{m}}K^{-\frac{1}{n}}_{1}K^{-\frac{2}{n}}_{2}\cdots K^{-\frac{n-1}{n}}_{n-1}K_{n+1}^{-\frac{m-1}{m}}\cdots
K_{n+m-1}^{-\frac{1}{m}}c^{-1}E_{n+m-1}\\
&=K_nE_{n+m-1}.
\end{split}
\end{equation*}

\emph{\textbf{Step 3}. The cross relations between the new simple-root vector $E_n$ and the $F_i$'s.}

The identity $[E_n,F_n]=\frac{K_n-K_n^{-1}}{q-q^{-1}}$ follows directly from $[e^{\underline{i}},f_{\underline{j}}]=
\frac{(m^{+})^{\underline{i}}_{\underline{j}}c^{-1}-c(m^{-})^{\underline{i}}_{\underline{}j}}{q_{\ast}-q_{\ast}^{-1}}$.

From the relations
\begin{gather*}
(m_{\mathbb{C}^n}^{-})^{i+1}_i=(q-q^{-1})(m_{\mathbb{C}^n}^{-})^{i+1}_{i+1}F_i, \qquad\qquad\quad 1\leq i\leq n-1,\\
(m_{(\mathbb{C}^m)^*}^{-})^{j+1}_j=(q-q^{-1})(m_{(\mathbb{C}^m)^*}^{-})^{j+1}_{j+1}F_{n+m-j},\quad 1\leq j\leq m-1.
\end{gather*}
we derive the following equalities:
\begin{equation}\label{A8}
(m^-)^{(i+1,m)}_{(i,m)}=(q-q^{-1})\,(m^-)^{(i+1,m)}_{(i+1,m)}\,F_i,\quad
(m^-)^{(n,j+1)}_{(n,j)}=(q-q^{-1})\,(m^-)^{(n,j+1)}_{(n,j+1)}\,F_{n+m-j}.
\end{equation}

This prompts us to focus on the cross relation $(m^{-})^{\underline{i}}_{\underline{j}}e^{\underline{k}}
=R_{VV}{}^{\underline{k}\,\underline{i}}_{\underline{a}\,\underline{b}}e^{\underline{a}}(m^{-})^{\underline{b}}_{\underline{j}}$ in the newly constructed quantum group, which yields the cross commutation relations between \(E_n\) and the remaining generators \(F_i\).
\begin{equation}\label{A9}
\left.
\begin{array}{rl}
(m^-)^{(i+1,m)}_{(i,m)}\,E_{n}&=(m^-)^{(i+1,m)}_{(i,m)}\,e^{(n,m)}
=(R_{\mathbb{C}^n\otimes(\mathbb{C}^m)^*})^{(nm)(i+1,m)}_{(ab)(cd)}\,e^{(a,b)}\,(m^-)^{(c,d)}_{(i,m)}\\
&=(R_{\mathbb{C}^n})^{n,i+1}_{ac}(R_{(\mathbb{C}^m)^*})^{mm}_{bd}\,e^{(a,b)}\,\Bigl((m_{\mathbb{C}^n}^-)^c_i\otimes(m_{(\mathbb{C}^m)^*}^-)^d_m\Bigr)\\
&=(R_{\mathbb{C}^n})^{n,i+1}_{n,i+1}(R_{(\mathbb{C}^m)^*})^{mm}_{mm}\,e^{(n,m)}\,\Bigl((m_{\mathbb{C}^n}^-)^{i+1}_i\otimes(m_{(\mathbb{C}^m)^*}^-)^m_m\Bigr)\\
&=(R_{\mathbb{C}^n})^{n,i+1}_{n,i+1}(R_{(\mathbb{C}^m)^*})^{mm}_{mm}\,E_{n}\,(m^-)^{(i+1,m)}_{(i,m)},
\end{array}
\right.
\end{equation}
\begin{equation}\label{A10}
\left.
\begin{array}{rl}
(m^-)^{(i+1,m)}_{(i+1,m)}\,E_{n}&=(m^-)^{(i+1,m)}_{(i+1,m)}\,e^{(n,m)}
=(R_{\mathbb{C}^n\otimes(\mathbb{C}^m)^*})^{(nm)(i+1,m)}_{(ab)(cd)}\,e^{(a,b)}\,(m^-)^{(c,d)}_{(i+1,m)}\\
&=(R_{\mathbb{C}^n})^{n,i+1}_{ac}(R_{(\mathbb{C}^m)^*})^{mm}_{bd}\,e^{(a,b)}\,\Bigl((m_{\mathbb{C}^n}^-)^c_{i+1}\otimes(m_{(\mathbb{C}^m)^*}^-)^d_m\Bigr)\\
&=(R_{\mathbb{C}^n})^{n,i+1}_{n,i+1}(R_{(\mathbb{C}^m)^*})^{mm}_{mm}\,e^{(n,m)}\,\Bigl((m_{\mathbb{C}^n}^-)^{i+1}_{i+1}\otimes(m_{(\mathbb{C}^m)^*}^-)^m_m\Bigr)\\
&=(R_{\mathbb{C}^n})^{n,i+1}_{n,i+1}(R_{(\mathbb{C}^m)^*})^{mm}_{mm}\,E_{n}\,(m^-)^{(i+1,m)}_{(i+1,m)}\,,
\end{array}
\right.
\end{equation}
Substituting the first equality of \eqref{A8} into equality \eqref{A9}, we obtain
\begin{equation}\label{A11}
\left.
\begin{array}{rl}
(m^-)^{(i+1,m)}_{(i+1,m)}F_iE_{n}=(R_{\mathbb{C}^n})^{n,i+1}_{n,i+1}(R_{(\mathbb{C}^m)^*})^{mm}_{mm}E_{n}(m^-)^{(i+1,m)}_{(i+1,m)}F_i
=\frac{(R_{\mathbb{C}^n})^{n,i+1}_{n,i+1}(R_{(\mathbb{C}^m)^*})^{mm}_{mm}}{(R_{\mathbb{C}^n})^{n,i+1}_{n,i+1}(R_{(\mathbb{C}^m)^*})^{mm}_{mm}}
(m^-)^{(i+1,m)}_{(i+1,m)}E_{n}F_i.
\end{array}
\right.
\end{equation}
Multiplying the both sides of equality \eqref{A11} on the left by $(m^+)^{(i+1,m)}_{(i+1,m)}$,
we prove that
\begin{equation*}
F_iE_{n}=E_{n}F_i,\quad
1\leq i\leq n-1.
\end{equation*}

On the other hand,
we also have
\begin{gather*}
(m^-)^{(n,i+1)}_{(n,i)}E_{n}=(m^-)^{(n,i+1)}_{(n,i)}e^{(n,m)}
=(R_{\mathbb{C}^n})^{nn}_{nn}(R_{(\mathbb{C}^m)^*})^{m,i+1}_{m,i+1}E_{n}(m^-)^{(n,i+1)}_{(n,i)},
\\
(m^-)^{(n,i+1)}_{(n,i+1)}E_{n}=(m^-)^{(n,i+1)}_{(n,i+1)}e^{(n,m)}
=(R_{\mathbb{C}^n})^{nn}_{nn}(R_{(\mathbb{C}^m)^*})^{m,i+1}_{m,i+1}E_{n}(m^-)^{(n,i+1)}_{(n,i+1)}.
\end{gather*}
Using the second equality  of \eqref{A8},
we get the cross relations
\begin{equation*}
F_{n+m-j}E_{n}=E_{n}F_{n+m-j},\quad
1\leq j\leq m-1.
\end{equation*}
In summary,
we prove that the following cross relations
\begin{equation}\label{A12}
[E_n,F_k]=\delta_{nk}\frac{K_n-K_n^{-1}}{q-q^{-1}},\quad \mbox{\it for~any~} k\in\{1,\cdots,n+m-1\}.
\end{equation}

\emph{\textbf{Step 4}. The $q$-Serre relations in the positive part.}

From the relation
$(m^+)^{i,m}_{i+1,m}=(q-q^{-1})E_i\,(m^+)^{i+1,m}_{i+1,m}$,
which is deduced from $(m_{\mathbb{C}^n}^+)^{i}_{i+1}=(q-q^{-1})E_i\,(m_{\mathbb{C}^n}^+)^{i+1}_{i+1}$,
we see that each $E_i$ occurs in the expression for $(m^+)^{i,m}_{i+1,m}$ for $1\leq i\leq n-1$.
Using the cross relation
\[
e^{\underline{i}}\,(m^{+})^{\underline{j}}_{\underline{k}}
=\lambda\, R^{\underline{j}\,\underline{i}}_{\underline{a}\,\underline{b}}\,(m^{+})^{\underline{a}}_{\underline{k}}\,e^{\underline{b}}
\]
in the newly constructed quantum group, we obtain the following equality:
\begin{equation}\label{A13}
\left.
\begin{array}{ll}
&E_n\,(m^+)^{i,m}_{i+1,m}
=e^{(n,m)}\,(m^+)^{i,m}_{i+1,m}
=\lambda R^{(im)(nm)}_{(ab)(cd)}\,(m^+)^{a,b}_{i+1,m}\,e^{(c,d)}\\
&
\quad=(R_{\mathbb{C}^n})^{in}_{ac}\,(R_{(\mathbb{C}^m)^*})^{mm}_{bd}\,\left((m_{\mathbb{C}^n}^+)^a_{i+1}\otimes (m_{(\mathbb{C}^m)^*}^+)^b_{m}\right)\,e^{(c,d)}
\\
&\quad=(R_{\mathbb{C}^n})^{in}_{in}\,(R_{(\mathbb{C}^m)^*})^{mm}_{mm}\,\left((m_{\mathbb{C}^n}^+)^i_{i{+}1}{\otimes}(m_{(\mathbb{C}^m)^*}^+)^m_{m}\right)\,e^{(n,m)}\\
&\qquad+\,
(R_{\mathbb{C}^n})^{in}_{i{+}1,n{-}1}\,(R_{(\mathbb{C}^m)^*})^{mm}_{mm}\,\left((m_{\mathbb{C}^n}^+)^{i{+}1}_{i{+}1}{\otimes}\,(m_{(\mathbb{C}^m)^*}^+)^m_{m}\right)\,e^{(n{-}1,m)}\\
&\quad=(R_{\mathbb{C}^n})^{in}_{in}\,(R_{(\mathbb{C}^m)^*})^{mm}_{mm}\,(m^+)^{i,m}_{i+1,m}\,e^{(n,m)}+
(R_{\mathbb{C}^n})^{in}_{i+1,n-1}\,(R_{(\mathbb{C}^m)^*})^{mm}_{mm}\,(m^+)^{i+1,m}_{i+1,m}\,e^{(n-1,m)}.
\end{array}
\right.
\end{equation}
Hence we obtain
\begin{gather}
e^{(n,m)}\,(m^+)^{i,m}_{i+1,m}=(R_{\mathbb{C}^n})^{in}_{in}\,(R_{(\mathbb{C}^m)^*})^{mm}_{mm}\,(m^+)^{i,m}_{i+1,m}\,e^{(n,m)}, \quad 1\leq i\leq n-2,\label{A14} \\
e^{(n,m)}(m^+)^{n-1,m}_{n,m}=(R_{\mathbb{C}^n})^{n-1,n}_{n-1,n}\,(R_{(\mathbb{C}^m)^*})^{mm}_{mm}\,(m^+)^{n-1,m}_{n,m}\,e^{(n,m)}{+}
(R_{\mathbb{C}^n})^{n-1,n}_{n,n-1}\,(R_{(\mathbb{C}^m)^*})^{mm}_{mm}\,(m^+)^{n,m}_{n,m}\,e^{(n-1,m)},  \label{A15}
\end{gather}
via the following entries of the associated $R$-matrix
$$
(R_{\mathbb{C}^n})^{i,n}_{i+1,n-1}=
q^{\delta_{in}}\delta_{i,i+1}\delta_{n,n-1}+(q-q^{-1})\,\delta_{i,n-1}\delta_{n,i+1}\theta(n-i)
=\left\{
\begin{array}{ll}
0, & i\neq n-1,\\
(q-q^{-1}), & i=n-1.
\end{array}
\right.
$$
We get the cross relation between $e^{(n,m)}$ and $(m^+)^{i+1,m}_{i+1,m}$ from equality \eqref{A4},
that is,
\begin{equation}\label{A16}
e^{(n,m)}(m^+)^{i+1,m}_{i+1,m}=\lambda R^{(i+1,m)(nm)}_{(ab)(cd)}\,(m^+)^{ab}_{i+1,m}\,e^{(c,d)}
=(R_{\mathbb{C}^n})^{i+1,n}_{i+1,n}\,(R_{(\mathbb{C}^m)^*})^{mm}_{mm}\,(m^+)^{i+1,m}_{i+1,m}\,e^{(n,m)}.
\end{equation}
Substituting $(m^+)^{i,m}_{i+1,m}=(q-q^{-1})\,E_i\,(m^+)^{i+1,m}_{i+1,m}$ into \eqref{A14}, we obtain
\begin{equation}\label{A17}
\left.
\begin{array}{rl}
E_n\,(q-q^{-1})\,E_i\,(m^+)^{i+1,m}_{i+1,m}&=(R_{\mathbb{C}^n})^{in}_{in}\,(R_{(\mathbb{C}^m)^*})^{mm}_{mm}\,(q-q^{-1})\,E_i\,(m^+)^{i+1,m}_{i+1,m}E_n
\\
&\overset{\eqref{A16}}=
\frac{(R_{\mathbb{C}^n})^{in}_{in}(R_{(\mathbb{C}^m)^*})^{mm}_{mm}}{(R_{\mathbb{C}^n})^{i+1,n}_{i+1,n}(R_{(\mathbb{C}^m)^*})^{mm}_{mm}}\,(q-q^{-1})\,E_i\,E_n\,(m^+)^{i+1,m}_{i+1,m}\\
&=(q-q^{-1})\,E_i\,E_n\,(m^+)^{i+1,m}_{i+1,m}.
\end{array}
\right.
\end{equation}
Namely,
we get the $q$-Serre relations:
\begin{equation}\label{A18}
E_nE_i=E_iE_n,\quad
1\leq i\leq n-2.
\end{equation}

According to $(m^+)^{n-1,m}_{n,m}=(q-q^{-1})\,E_{n-1}\,(m^+)^{n,m}_{n,m}$ and \eqref{A15}, we obtain
\begin{equation}\label{A19}
\left.
\begin{array}{rl}
&e^{(n,m)}\,(q-q^{-1})\,E_{n-1}\,(m^+)^{n,m}_{n,m}\\
=&(R_{\mathbb{C}^n})^{n-1,n}_{n-1,n}(R_{(\mathbb{C}^m)^*})^{mm}_{mm}\,(q-q^{-1})\,E_{n-1}\,(m^+)^{n,m}_{n,m}\,e^{(n,m)}{+}
(R_{\mathbb{C}^n})^{n{-}1,n}_{n,n{-}1}(R_{(\mathbb{C}^m)^*})^{mm}_{mm}\,(m^+)^{n,m}_{n,m}\,e^{(n{-}1,m)}\\
=&1\cdot q\cdot(q-q^{-1})\,E_{n-1}\,(m^+)^{n,m}_{n,m}\,e^{(n,m)}+(q-q^{-1})\cdot q\,(m^+)^{n,m}_{n,m}\,e^{(n-1,m)}.
\end{array}
\right.
\end{equation}
Then the equality
\begin{equation}\label{A20}
e^{(n,m)}E_{n-1}(m^+)^{n,m}_{n,m}=qq^{-2}E_{n-1}\,e^{(n,m)}\,(m^+)^{n,m}_{n,m}+qq^{-1}\,e^{(n-1,m)}\,(m^+)^{n,m}_{n,m}
\end{equation}
follows from \eqref{A19}, by virtue of the cross relations
\begin{gather*}
e^{(n-1,m)}\,(m^+)^{n,m}_{n,m}=\lambda R^{(nm)(n-1,m)}_{(ab)(cd)}\,(m^+)^{a,b}_{n,m}\,e^{(c,d)}=q\,(m^+)^{n,m}_{n,m}\,e^{(n-1,m)},
\\
e^{(n,m)}\,(m^+)^{n,m}_{n,m}=\lambda R^{(nm)(nm)}_{(ab)(cd)}\,(m^+)^{a,b}_{n,m}\,e^{(c,d)}=q^2\,(m^+)^{n,m}_{n,m}\,e^{(n,m)}.
\end{gather*}
Consequently, we obtain
\begin{equation}\label{A21}
e^{(n-1,m)}=e^{(n,m)}\,E_{n-1}-q^{-1}\,E_{n-1}\,e^{(n,m)}.
\end{equation}

To obtain the relation between $E_n$ and $E_{n-1}$, it suffices to determine the relations among $e^{(n-1,m)}$, $e^{(n,m)}$ and $E_{n-1}$.
Based on the fact that $E_{n-1}$ locates at the entry $(m^+)^{n-1,m}_{n,m}$, we consider the equality
\begin{equation}\label{A22}
e^{(n-1,m)}\,(m^+)^{n-1,m}_{n,m}=\lambda R^{(n-1,m)(n-1,m)}_{(ab)(cd)}\,(m^+)^{a,b}_{n,m}\,e^{(c,d)}=q\,(m^+)^{n-1,m}_{n,m}\,e^{(n-1,m)}.
\end{equation}
Then, using $(m^+)^{n-1,m}_{n,m}=(q-q^{-1})E_{n-1}(m^+)^{n,m}_{n,m}$ and $e^{(n-1,m)}(m^+)^{n,m}_{n,m}=q(m^+)^{n,m}_{n,m}e^{(n-1,m)}$, we derive the following cross relation \eqref{A23} from \eqref{A22}:
\begin{equation}\label{A23}
e^{(n-1,m)}\,E_{n-1}=qE_{n-1}\,e^{(n-1,m)}.
\end{equation}

Substituting this equality into \eqref{A21}, we obtain the desired $q$-Serre relation
\begin{equation}\label{A24}
E_nE_{n-1}^2-(q+q^{-1})E_{n-1}E_nE_{n-1}+E_{n-1}^2E_n=0.
\end{equation}

On the other hand, from \eqref{rmatrix1} we get the non-vanishing entries in the row $(n,m)(n-1,m)$ of the matrix $R'$, namely
\[
R'^{(n,m)(n-1,m)}_{(n-1,m)(n,m)}=q^2+2,\qquad
R'^{(n,m)(n-1,m)}_{(n,m)(n-1,m)}=-(q+q^{-1}).
\]
This yields
\[
e^{(n-1,m)}\,e^{(n,m)}=R'^{(nm)(n-1,m)}_{(ab)(cd)}\,e^{(a,b)}\,e^{(c,d)}
=(q^2+2)\,e^{(n-1,m)}\,e^{(n,m)}-(q+q^{-1})\,e^{(n,m)}\,e^{(n-1,m)},
\]
which further gives $e^{(n,m)}\,e^{(n-1,m)}=q\,e^{(n-1,m)}\,e^{(n,m)}$.
Applying \eqref{A21} again, we obtain another $q$-Serre relation
\begin{equation}\label{A25}
E_n^2E_{n-1}-(q+q^{-1})E_nE_{n-1}E_n+E_{n-1}E_n^2=0.
\end{equation}

So far, all $q$-Serre relations for the positive part have been proved.
The corresponding $q$-Serre relations for the newly added simple root vector $F_n$ and the remaining $F_i$'s can be verified by analogous arguments.
We thus prove that the algebraic structure is compatible with the identification given in Theorem \ref{typeA}.
Compatibility of the comultiplication and antipode under this identification can be checked easily.
Hence, the identification in the Theorem is a Hopf algebra morphism.

\emph{\textbf{Step 5}. A vector space basis of braided group $V(R',R)$.}

Denote by $\Delta=\{\alpha_1,\cdots,\alpha_{n+m-1}\}$ the simple roots of $\mathfrak{sl}_{n+m}$,
$\Phi$ its the root system,
then the subset $\textrm{I}=\{\alpha_1,\cdots,\alpha_{n-1},\alpha_{n+1},\cdots,\alpha_{n+m-1}\}\subset\Delta$ defines a root subsystem $\Phi_{\textrm{I}}\subset\Phi$,
with positive roots $\Phi^{+}_{\textrm{I}}$,
negative roots $\Phi^{-}_{\textrm{I}}$.
As is well known,
$\mathfrak{sl}_{n+m}=\mathfrak{u}^{-}_\textrm{I}\oplus(\mathfrak{h}\oplus\sum\limits_{\alpha\in\Phi_{\textrm{I}}}\mathfrak{g}_{\alpha})\oplus\mathfrak{u}_\textrm{I}$,
and $\mathfrak{u}^{-}_\textrm{I}=\bigoplus_{\alpha\in\Phi^{-}\setminus\Phi^{-}_{\textrm{I}}}\mathfrak{g}_{\alpha}$,
$\mathfrak{u}_\textrm{I}=\bigoplus_{\alpha\in\Phi^{+}\setminus\Phi^{+}_{\textrm{I}}}\mathfrak{g}_{\alpha}$.

With $|\Phi^{-}\setminus\Phi^{-}_{\textrm{I}}|=|\Phi^{+}\setminus\Phi^{+}_{\textrm{I}}|=nm$ and PBW-basis of $U_q(\mathfrak{sl}_{n+m})$,
we need to prove that $V(R',R)\cong S(\mathbb{C}^{nm})$ (the symmetric algebra) as vector spaces,
to make sure that the identification in Theorem \ref{typeA} is a bijection.

To this end, it suffices to construct a vector-space basis for $V(R',R)$.
We then invoke the Diamond lemma from ring theory (\cite{Ber}).
The algebra $V(R',R)$ has $nm$ generators $e^{\underline{i}}$,
where $\underline{i}=(i_1,i_2)$ with
$i_1\in \textrm{I}_1=\{1,\dots,n\}$ and
$i_2\in\textrm{I}_2=\{1,\dots,m\}$.
We introduce the following ordering on these indices $\underline{i}$ as follows:
\begin{gather*}
\underline{i}=(i_1,i_2)>\underline{j}=(j_1,j_2),\mbox{\quad if~} i_1=j_1, i_2>j_2 \mbox{~or~} i_1>j_1,\\
\underline{i}=(i_1,i_2)<\underline{j}=(j_1,j_2),\mbox{\quad if~} i_1=j_1, i_2<j_2\mbox{~or~}i_1<j_1.
\end{gather*}

According to $R_{V}{}^{ii}_{ii}=q, R_{V}{}^{ij}_{ij}=1$,
and $R_{V}{}^{ij}_{ji}=q-q^{-1}$ if $i<j$,
where $V$ denotes $\mathbb{C}^n$ or $(\mathbb{C}^{m})^*$,
we get
$$
(R_{V}PR_{V})^{ii}_{ii}=q^2,\quad
(R_{V}PR_{V})^{ij}_{ij}=q-q^{-1},\quad
(R_{V}PR_{V})^{ij}_{ji}= \left\{
\begin{array}{lcl}
	1+(q-q^{-1})^2&~~&i<j.\\
	1&~~&i>j,
\end{array}
\right.
$$
This forces us to examine the following cases for the relation
$e^{\underline{j}}e^{\underline{i}}=\sum_{\underline{a},\underline{b}} R'{}^{\underline{i}\underline{j}}_{\underline{a}\,\underline{b}}e^{\underline{a}}e^{\underline{b}}$.

\emph{\textbf{Case 1}: $\underline{i}>\underline{j}$ with  $i_1=j_1, i_2>j_2$, then $\underline{i}=(i_1,i_2)$ and $\underline{j}=(i_1,j_2)$.}
$$
\left.
\begin{array}{rl}
R'{}^{\underline{i}\underline{j}}_{\underline{a}\,\underline{b}}=&
(R_{\mathbb{C}^n}PR_{\mathbb{C}^n})^{i_1i_1}_{a_1b_1}(R_{(\mathbb{C}^m)^*}PR_{(\mathbb{C}^m)^*})^{i_2j_2}_{a_2b_2}
-(q^2+q^{-2})(R_{\mathbb{C}^n})^{i_1i_1}_{a_1b_1}(R_{(\mathbb{C}^m)^*})^{i_2j_2}_{a_2b_2}+2P^{\underline{i}\underline{j}}_{\underline{a}\,\underline{b}}\\
=&
\left\{\begin{array}{ll}
R'{}^{(i_1i_1)(i_1j_2)}_{(i_1i_1)(i_1j_2)}=-(q+q^{-1}),& \mbox{\it when}\quad (a_1b_1)=(i_1i_1), \ (a_2b_2)=(i_2j_2); \\
R'{}^{(i_1i_1)(i_1j_2)}_{(i_1j_2)(i_1i_2)}=q^2+2,&\mbox{\it when}\quad (a_1b_1)=(i_1i_1), \ (a_2b_2)=(j_2i_2).
  \end{array}
\right.
\end{array}
\right.
$$
With these results in hand,
we obtain 
$$
e^{(i_1j_2)}e^{(i_1i_2)}=-(q+q^{-1})e^{(i_1i_2)}e^{(i_1j_2)}+(q^2+2)e^{(i_1j_2)}e^{(i_1i_2)},
$$
which yields the following equality in this case:
\begin{equation}\label{alg0}
e^{(i_1i_2)}e^{(i_1j_2)}=qe^{(i_1j_2)}e^{(i_1i_2)},\quad\mbox{i.e.},\quad e^{\underline{i}}e^{\underline{j}}=qe^{\underline{j}}e^{\underline{i}}.
\end{equation}

\emph{\textbf{Case 2}: $\underline{i}>\underline{j}$ with  $i_1>j_1, i_2=j_2$, or $i_2>j_2$, or $i_2<j_2$}.
$$R'{}^{\underline{i}\underline{j}}_{\underline{a}\,\underline{b}}=
(R_{\mathbb{C}^n}PR_{\mathbb{C}^n})^{i_1j_1}_{a_1b_1}(R_{(\mathbb{C}^m)^*}PR_{(\mathbb{C}^m)^*})^{i_2j_2}_{a_2b_2}
-(q^2+q^{-2})(R_{\mathbb{C}^n})^{i_1j_1}_{a_1b_1}(R_{(\mathbb{C}^m)^*})^{i_2j_2}_{a_2b_2}+2P^{\underline{i}\underline{j}}_{\underline{a}\,\underline{b}}.$$
Thus we have
\begin{itemize}
\setlength\itemsep{0em}
\item [{(1)}] \emph{when $(a_1b_1)=(i_1j_1)$, $(a_2b_2)=(i_2j_2)$,}
$$
\left.
\begin{array}{rl}
R'{}^{(i_1i_2)(j_1j_2)}_{(i_1i_2)(j_1j_2)}
=&(R_{\mathbb{C}^n}PR_{\mathbb{C}^n})^{i_1j_1}_{i_1j_1}(R_{(\mathbb{C}^m)^*}PR_{(\mathbb{C}^m)^*})^{i_2j_2}_{i_2j_2}
-(q^2+q^{-2})(R_{\mathbb{C}^n})^{i_1j_1}_{i_1j_1}(R_{(\mathbb{C}^m)^*})^{i_2j_2}_{i_2j_2}\\
=&(q-q^{-1})(R_{(\mathbb{C}^m)^*}PR_{(\mathbb{C}^m)^*})^{i_2j_2}_{i_2j_2}-(q^2+q^{-2})(R_{(\mathbb{C}^m)^*})^{i_2j_2}_{i_2j_2}\\
=&
\left\{
\begin{array}{ll}
(q-q^{-1})(q-q^{-1})-(q^2+q^{-2})=-2,&\mbox{if~} i_2\neq j_2,\\
(q-q^{-1})q^2-(q^2+q^{-2})q=-(q+q^{-1}),&\mbox{if~} i_2=j_2.
\end{array}
\right.
\end{array}
\right.
$$
\item [{(2)}]\emph{ when $(a_1b_1)=(i_1j_1)$, $(a_2b_2)=(j_2i_2)$,}
$$
\left.
\begin{array}{rl}
R'{}^{(i_1i_2)(j_1j_2)}_{(i_1j_2)(j_1i_2)}
=&(R_{\mathbb{C}^n}PR_{\mathbb{C}^n})^{i_1j_1}_{i_1j_1}(R_{(\mathbb{C}^m)^*}PR_{(\mathbb{C}^m)^*})^{i_2j_2}_{j_2i_2}
-(q^2+q^{-2})(R_{\mathbb{C}^n})^{i_1j_1}_{i_1j_1}(R_{(\mathbb{C}^m)^*})^{i_2j_2}_{j_2i_2}\\
=&\left\{
\begin{array}{ll}
(q-q^{-1}).1-(q^2+q^{-2}).0=(q-q^{-1}),&\mbox{if~} i_2>j_2,\\
(q-q^{-1})[1+(q-q^{-1})^2]-(q^2+q^{-2})(q-q^{-1})=-(q-q^{-1}),&\mbox{if~} i_2<j_2,\\
(q-q^{-1})q^2-(q^2+q^{-2})q=-(q+q^{-1}),&\mbox{if~} i_2=j_2.
\end{array}
\right.
\end{array}
\right.
$$
\item [{(3)}] \emph{when $(a_1b_1)=(j_1i_1)$, $(a_2b_2)=(i_2j_2)$.}
In this case,
$\underline{a}=(j_1,i_2)$, $\underline{b}=(i_1,j_2)$.
Then $\underline{a}=\underline{j}, \underline{b}=\underline{i}$ whenever $i_2=j_2$; thus we have
$$\left.
\begin{array}{rl}
R'{}^{(i_1i_2)(j_1j_2)}_{(j_1i_2)(i_1j_2)}{=}&
(R_{\mathbb{C}^n}PR_{\mathbb{C}^n})^{i_1j_1}_{j_1i_1}(R_{(\mathbb{C}^m)^*}PR_{(\mathbb{C}^m)^*})^{i_2j_2}_{i_2j_2}
{-}(q^2+q^{-2})(R_{\mathbb{C}^n})^{i_1j_1}_{j_1i_1}(R_{(\mathbb{C}^m)^*})^{i_2j_2}_{i_2j_2}+2P^{\underline{i}\underline{j}}_{\underline{a}\,\underline{b}}\\
{=}&(R_{(\mathbb{C}^m)^*}PR_{(\mathbb{C}^m)^*})^{i_2j_2}_{i_2j_2}+2P^{\underline{i}\underline{j}}_{\underline{a}\underline{b}}\\
{=}&\left\{
\begin{array}{ll}
q-q^{-1},&\mbox{if~} i_2\neq j_2,\\
q^2+2,&\mbox{if~} i_2= j_2.
\end{array}
\right.
\end{array}
\right.
$$
\item [{(4)}]\emph{ when $(a_1b_1)=(j_1i_1)$, $(a_2b_2)=(j_2i_2)$ ,}
$$
\left.
\begin{array}{rl}
R'{}^{(i_1i_2)(j_1j_2)}_{(j_1j_2)(i_1i_2)}
=&(R_{\mathbb{C}^n}PR_{\mathbb{C}^n})^{i_1j_1}_{j_1i_1}(R_{(\mathbb{C}^m)^*}PR_{(\mathbb{C}^m)^*})^{i_2j_2}_{j_2i_2}
-(q^2+q^{-2})(R_{\mathbb{C}^n}{}^{i_1j_1}_{j_1i_1})(R_{(\mathbb{C}^m)^*})^{i_2j_2}_{j_2i_2}+2\\
=&(R_{(\mathbb{C}^m)^*}PR_{(\mathbb{C}^m)^*})^{i_2j_2}_{j_2i_2}+2=
\left\{
\begin{array}{ll}
3,&\mbox{if~} i_2>j_2,\\
3+(q-q^{-1})^2,&\mbox{if~} i_2<j_2,\\
q^2+2,&\mbox{if~} i_2=j_2.
\end{array}
\right.\\
\end{array}
\right.
$$
\end{itemize}

With these results in hand,
we obtain the following relations in the braided group $V(R',R)$:
\begin{equation}\label{B0}
e^{(i_1,i_2)}e^{(j_1,j_2)}=\frac{q-q^{-1}}{2}(e^{(i_1,j_2)}e^{(j_1,i_2)}+e^{(j_1,i_2)}e^{(i_1,j_2)})+e^{(j_1,j_2)}e^{(i_1,i_2)},\quad\mbox{\it for}\ i_2>j_2,
\end{equation}
\begin{equation}\label{B1}
e^{(i_1,i_2)}e^{(j_1,j_2)}{=}\frac{q{-}q^{-1}}{2}(e^{(j_1,i_2)}e^{(i_1,j_2)}{-}e^{(i_1,j_2)}e^{(j_1,i_2)}){+}\frac{2{+}
(q{-}q^{-1})^2}{2}e^{(j_1,j_2)}e^{(i_1,i_2)},\quad\mbox{\it for}\ i_2<j_2,
\end{equation}
\begin{equation}\label{B2}
e^{\underline{i}}e^{\underline{j}}=\frac{2q^2+3}{2(q+q^{-1})}e^{\underline{j}}e^{\underline{i}},\quad\mbox{\it for}\ i_2=j_2.
\end{equation}

Using relation \eqref{B1}, we derive the following relation from \eqref{B0}:
\begin{equation}\label{B3}
e^{\underline{i}}e^{\underline{j}}=\frac{4+(q-q^{-1})^2}{4+2(q+q^{-1})}\left[e^{\underline{j}}e^{\underline{i}}
+(q-q^{-1})e^{(j_1,i_2)}e^{(i_1,j_2)}\right],\quad\mbox{\it for}\ i_2>j_2.
\end{equation}
Moreover,
$i_1>j_1$ means $(j_1,i_2)<(i_1,j_2)$,
and $i_2>j_2$ means $(j_1,i_2)<\underline{i}$ and $(i_1,j_2)<\underline{i}$.
Hence, by relations \eqref{alg0}, \eqref{B2} and \eqref{B3} in the braided group $V(R',R)$, the set of monomials $\{\,(e^{1,1})^{k_{11}}\cdots(e^{1,m})^{k_{1m}}(e^{2,1})^{k_{21}}\cdots(e^{2,m})^{k_{2m}}\cdots(e^{n,1})^{k_{n1}}\cdots(e^{n,m})^{k_{nm}}\mid k_{ij}=0,1,2,\cdots\}$ forms a vector‑space basis of $V(R',R)$.
Thus $V(R',R)$ is isomorphic to $S(\mathbb{C}^{nm})$ as a vector space.

Together with the foregoing arguments, this establishes that the identification in the theorem is a Hopf algebra isomorphism.
\end{proof}

\subsection{Non-simply-laced case: Type $F_4$}
We will give the grafting construction of the quantum enveloping algebra $U_q(F_4)$ as a sample in the non-simply-laced case in this subsection.

To do this, we choose $\mathfrak{g}_{1}\simeq\mathfrak{sl}_{3}$ with simple roots $\alpha_i$ and $(\alpha_i,\alpha_i)=2,\,i=1,2$,
and $\mathfrak{g}_{2}\simeq\mathfrak{sl}_{2}$ with the simple root $\alpha_4$ and $(\alpha_4,\alpha_4)=1$.
Denote by $\{E_i, F_i, K_i\mid i=1,2\}$ the generators of the quantum enveloping algebra $U_q(\mathfrak{g}_1)$,
and $\{E_4, F_4, K_4\}$ the generators of $U_q(\mathfrak{g}_2)$.
We choose the $3$-dimensional natural representation $\mathbb{C}^3$ of $U_q(\mathfrak{g}_1)$ and the $2$-dimensional natural representation $\mathbb{C}^2$ of $U_q(\mathfrak{g}_2)$.

The corresponding $R$-matrices are as follows.
$$
R_{\mathbb{C}^3}=
\left(
\begin{array}{ccccccccc}
q&~0&~0&~0&~0&~0&~0&~0&~0\\
0&~1&~0&~q-q^{-1}&~0&~0&~0&~0&~0\\
0&~0&~1&~0&~0&~0&~q-q^{-1}&~0&~0\\
0&~0&~0&~1&~0&~0&~0&~0&~0\\
0&~0&~0&~0&~q&~0&~0&~0&~0\\
0&~0&~0&~0&~0&~1&~0&~q-q^{-1}&~0\\
0&~0&~0&~0&~0&~0&~1&~0&~0\\
0&~0&~0&~0&~0&~0&~0&~1&~0\\
0&~0&~0&~0&~0&~0&~0&~0&~q
\end{array}
\right),\qquad
R_{\mathbb{C}^2}=
\left(
\begin{array}{cccc}
q^{\frac{1}{2}}&~0&~0&~0\\
0&~1&~q^{\frac{1}{2}}-q^{-\frac{1}{2}}&~0\\
0&~0&~1&~0\\
0&~0&~0&~q^{\frac{1}{2}}
\end{array}
\right).
$$
Then, by Proposition~\ref{minimal}, the braiding $PR_{\mathbb{C}^3\otimes\mathbb{C}^2}$ associated with the representation $\mathbb{C}^3\otimes\mathbb{C}^2$ of $U_q(\mathfrak{g}_{1})\otimes U_q(\mathfrak{g}_{2})$ is diagonalizable, with distinct eigenvalues $q^{\frac32}$, $-q^{\frac12}$, $-q^{-\frac12}$, and $q^{-\frac32}$.
Hence its minimal polynomial equation is
\begin{equation*}
(PR_{\mathbb{C}^3\otimes\mathbb{C}^2}-q^{\frac{3}{2}}\textrm{I})(PR_{\mathbb{C}^3\otimes\mathbb{C}^2}+q^{\frac{1}{2}}\textrm{I})
(PR_{\mathbb{C}^3\otimes\mathbb{C}^2}+q^{-\frac{1}{2}}\textrm{I})(PR_{\mathbb{C}^3\otimes\mathbb{C}^2}-q^{-\frac{3}{2}}\textrm{I})=0.
\end{equation*}
To construct $U_q(F_4)$,
we take the Majid pair $(R, R')$ as given by
\begin{equation}\label{rmatrixF4}
R=q^{-\frac{1}{2}}R_{\mathbb{C}^3\otimes\mathbb{C}^2},\quad
R'=R(PR)^2+(q^{-1}-q-q^{-2})R(PR)+(q^{-1}-1-q^{-3})R+(q^{-2}+1)P.
\end{equation}

Also,
we obtain the dually-paired braided groups $V^{\vee}(R^{\prime},R_{21}^{-1})$ and $V(R',R)$.
Moreover,
we obtain the following lemma, whose proof is similar to that of Proposition 3.1 in \cite{HH3}.
\begin{lemma}\label{radi}
For the dual pair $V^{\vee}(R^{\prime},R_{21}^{-1})$ and $V(R',R)$,
we have
\begin{gather*}
\langle f,\,\left(e^{(i+1,j)}\right)^2e^{(i,j)}-(1+q)e^{(i+1,j)}e^{(i,j)}e^{(i+1,j)}+qe^{(i,j)}\left(e^{(i+1,j)}\right)^2\rangle=0,
 \\
\langle f,\,\left(e^{(3,j)}\right)^2e^{(1,j)}-(1+q)e^{(3,j)}e^{(1,j)}e^{(3,j)}+qe^{(1,j)}\left(e^{(3,j)}\right)^2\rangle=0,\quad
\mbox{for all} ~f\in V^{\vee}(R^{\prime},R_{21}^{-1}).\\
\langle f_{(i,j)}\left(f_{(i+1,j)}\right)^{2}-(1+q^{-1})f_{(i+1,j)}f_{(i,j)}f_{(i+1,j)}+q^{-1}\left(f_{(i+1,j)}\right)^{2}f_{(i,j)},\,e\rangle=0,\\
\langle f_{(1,j)}\left(f_{(3,j)}\right)^{2}-(1+q^{-1})f_{(3,j)}f_{(1,j)}f_{(3,j)}+q^{-1}\left(f_{(3,j)}\right)^{2}f_{(1,j)},\,e\rangle=0,\quad
\mbox{for all} ~e\in V(R',R).
\end{gather*}
\end{lemma}
Define $\tilde{V}(R',R), \tilde{V}^{\vee}(R^{\prime},R_{21}^{-1})$ as the quotient algebras of $V(R',R),V^{\vee}(R^{\prime},R_{21}^{-1})$ modulo the ideals generated by the above radicals, respectively.
Via some complicated calculations, we are able to obtain the following structural-relations diagram for $\tilde{V}(R',R)$, as well as
the counterpart for $\tilde{V}^{\vee}(R^{\prime},R_{21}^{-1})$ as described below.
\begin{center}
	\setlength{\unitlength}{1mm}
	\begin{picture}(55,75)
		\put(9,4){\circle*{1}}
		\put(2,20){\mbox{$(\blacktriangle$)}}
		\put(2,50){\mbox{$(\blacktriangle$)}}
		\put(42,20){\mbox{$(\blacktriangle$)}}
		\put(42,50){\mbox{$(\blacktriangle$)}}
		\put(1,3){$e^{(3,1)}$}
		\put(10,5){\line(1,1){30}}
		\put(18,5){\mbox{q-commu}}
		\put(10,4){\line(1,0){30}}
		\put(41,4){\circle*{1}}
		\put(43,3){$e^{(3,2)}$}
		\put(9,5){\line(0,1){30}}
		\put(41,5){\line(0,1){30}}
		\put(9,36){\circle*{1}}
		\put(1,36){$e^{(2,1)}$}
		\put(10,35){\line(1,-1){30}}
		\put(10,68){\line(1,-2){30.5}}
		\put(10,37){\line(1,1){30}}
		\put(10,5){\line(1,2){30.5}}
		\put(18,37.5){\mbox{q-commu}}
		\put(24.5,35){$\blacksquare$}
		\put(10,36){\line(1,0){30}}
		\put(41,36){\circle*{1}}
		\put(43,36){$e^{(2,2)}$}
		\put(10,68){\line(1,-1){30}}
		\put(24.5,51.5){$\blacksquare$}
		\put(24,19){$\blacksquare$}
		\put(9,37){\line(0,1){30}}
		\put(41,37){\line(0,1){30}}
		\put(9,69){\circle*{1}}
		\put(1,69){$e^{(1,1)}$}
		\put(18,70){\mbox{q-commu}}
		\put(10,69){\line(1,0){30}}
		\put(41,69){\circle*{1}}
		\put(43,69){$e^{(1,2)}$}
	\end{picture}
\end{center}
The relations along the horizontal lines are: $$e^{(j,2)}e^{(j,1)}=q^{\frac{1}{2}}e^{(j,1)}e^{(j,2)}, \quad j=1,2,3.$$ The relations along the three diagonal lines are:
$$
e^{(i,1)}e^{(i+1,2)}-q^{-\frac{1}{2}}e^{(i,2)}e^{(i+1,1)}+qe^{(i+1,1)}e^{(i,2)}-q^{\frac{1}{2}}e^{(i+1,2)}e^{(i,1)}=0, \quad i=1,2.
\leqno{(\blacksquare)}$$
The relations among the generators on the vertical lines are:
$$
\left\{
\begin{array}{l}
\left(e^{(i+1,j)}\right)^2e^{(i,j)}+qe^{(i,j)}\left(e^{(i+1,j)}\right)^2=(1+q)e^{(i+1,j)}e^{(i,j)}e^{(i+1,j)},\quad i=1,2,j=1,2,\\
\left(e^{(3,j)}\right)^2e^{(1,j)}+qe^{(1,j)}\left(e^{(3,j)}\right)^2=(1+q)e^{(3,j)}e^{(1,j)}e^{(3,j)},\quad j=1,2.
\end{array}
\right.
\leqno{(\blacktriangle)}
$$

Based on these,
we obtain the graft of quantum enveloping algebra $U_q(F_{4})$ on the new quantum group $U(\tilde{V}^{\vee}(R^{\prime},R_{21}^{-1}),U_q^{\text{ext}}(\mathfrak{g}_1\oplus\mathfrak{g}_2)\otimes
k[c,c^{-1}],\tilde{V}(R^{\prime},R))$ by Theorem \ref{multi}.
\begin{theorem}\label{typeF}
Upon identifying the braided-group generators $e^{(3,2)}$, $f_{(3,2)}$ and the group-like element $(m^+)^{(3,2)}_{(3,2)}c^{-1}$ with the additional simple root vectors $E_3$, $F_3$, $K_3$, respectively,
the resulting quantum group
\[
U=U\Bigl(\tilde{V}^{\vee}(R^{\prime},R_{21}^{-1}),U_q^{\text{ext}}(\mathfrak{sl}_3\oplus\mathfrak{sl}_2)\otimes k[c,c^{-1}],\tilde{V}(R^{\prime},R)\Bigr)
\]
is isomorphic to the quantum enveloping algebra $U_q(F_{4})$ with $K_i^{\pm\frac12}$ ($i=1,2$), $K_4^{\frac12}$ adjoined.
The extra simple root is represented by the filled circle in the Dynkin diagram below:
\begin{center}
\setlength{\unitlength}{1mm}
\begin{picture}(40,10)
\put(2,6){\circle{1}}
\put(0,8){$\alpha_1$}
\put(2.5,6){\line(1,0){10}}
\put(13,6){\circle{1}}
\put(12,8){$\alpha_2$}
\multiput(13.5,5.7)(1,0){10}{\line(1,0){0.5}}
\multiput(13.5,6.3)(1,0){10}{\line(1,0){0.5}}
\put(16.5,5){$>$}
\put(23.5,6){\circle*{1}}
\multiput(23.5,6)(1,0){10}{\line(1,0){0.5}}
\put(5.5,8){$\mathfrak{sl}_{3}$}
\put(25.5,8){$\mathfrak{sl}_{2}$}
\put(33.5,6){\circle{1}}
\put(4,0){Figure 3: Type $F_4$}
\put(32,8){$\alpha_{4}$}
\end{picture}
\end{center}
\end{theorem}
\begin{proof}
The proof is similar to that of Theorem \ref{typeA}. Here we give mainly the relations in the positive part of the new quantum group $U$ and a basis of the braided group $\tilde{V}(R',R)$.
The matrices $m^{\pm}$ consisting of the FRT-generators for the natural representation $\mathbb{C}^3$ of the quantum enveloping algebra $U_q(\mathfrak{sl}_3)$ are as follows.
$$
\left.
\begin{array}{rcl}
m_{\mathbb{C}^3}^{+}=
\left(
\begin{array}{ccc}
K^{\frac{2}{3}}_{1}K^{\frac{1}{3}}_{2}
&~~(q-q^{-1})
E_{1}K^{-\frac{1}{3}}_{1}K^{\frac{1}{3}}_{2}
&~~q^{-1}(q-q^{-1})
[E_1, E_2]_{q^{-1}}K^{-\frac{1}{3}}_{1}K^{-\frac{2}{3}}_{2}\\
0&K^{-\frac{1}{3}}_{1}K^{\frac{1}{3}}_{2}
&(q-q^{-1})
E_{2}K^{-\frac{1}{3}}_{1}K^{-\frac{2}{3}}_{2}\\
0&0&K^{-\frac{1}{3}}_{1}K^{-\frac{2}{3}}_{2}
\end{array}
\right)
\end{array}
\right.,
$$
$$
\left.
\begin{array}{rcl}
m_{\mathbb{C}^3}^{-}=
\left(
\begin{array}{ccc}
K^{-\frac{2}{3}}_{1}K^{-\frac{1}{3}}_{2}&0&0\\
(q-q^{-1})
K^{\frac{1}{3}}_{1}K^{-\frac{1}{3}}_{2}F_{1}
&K^{\frac{1}{3}}_{1}K^{-\frac{1}{3}}_{2}&0\\
q(q-q^{-1})
K^{\frac{1}{3}}_{1}K^{\frac{2}{3}}_{2}[F_2, F_1]_q
&~~(q-q^{-1})
K^{\frac{1}{3}}_{1}K^{\frac{2}{3}}_{2}F_{2}
&~~K^{\frac{1}{3}}_{1}K^{\frac{2}{3}}_{2}
\end{array}
\right)
\end{array}
\right.,
$$
where $[E_1, E_2]_{q^{-1}}=E_1E_2-q^{-1}E_2E_1$, $[F_2, F_1]_q=F_2F_1-qF_1F_2$.
The matrices $m^{\pm}$ for the natural representation $\mathbb{C}^2$ of quantum enveloping algebras $U_q(\mathfrak g_2)$ are:
$$
m_{(\mathbb{C}^2)^*}^{+}=
\left(
\begin{array}{cc}
K_4^{\frac{1}{2}}
&~~(q^{\frac{1}{2}}-q^{-\frac{1}{2}})
E_{4}K_4^{-\frac{1}{2}}\\
0&K_4^{-\frac{1}{2}}
\end{array}
\right), \qquad
m_{(\mathbb{C}^2)^*}^{-}=
\left(
\begin{array}{cc}
K_4^{-\frac{1}{2}}&0\\
(q^{\frac{1}{2}}-q^{-\frac{1}{2}})K_4^{\frac{1}{2}}F_4&K_4^{\frac{1}{2}}
\end{array}
\right).
$$

First, the main cross relations in the positive part of the new quantum group $U$ are as follows:
\begin{enumerate} \item The cross relations between the new generator $e^{(3,2)}$ and $K_i,i=1,2,3,4$.
$$
\left.
\begin{array}{rl}
e^{(3,2)}(m^+)^{(3,2)}_{(3,2)}c^{-1}=&\lambda R^{(3,2)(3,2)}_{(a_1,a_2)(b_1,b_2)}(m^+)^{(a_1,a_2)}_{(3,2)}e^{(b_1,b_2)}c^{-1}\\
=&\lambda R^{(3,2)(3,2)}_{(a_1,a_2)(b_1,b_2)}(m^+)^{(a_1,a_2)}_{(3,2)}\frac{1}{\lambda}c^{-1}e^{(b_1,b_2)}\\
=&q^{-\frac{1}{2}}R_{\mathbb{C}^3}{}^{(3,3)}_{(3,3)}R_{(\mathbb{C}^2)^*}{}^{(2,2)}_{(2,2)}(m^+)^{(3,2)}_{(3,2)}c^{-1}e^{(3,2)}\\
=&q^{-\frac{1}{2}}qq^{\frac{1}{2}}(m^+)^{(3,2)}_{(3,2)}c^{-1}e^{(3,2)}=q(m^+)^{(3,2)}_{(3,2)}c^{-1}e^{(3,2)}.
\end{array}
\right.
$$
On the one hand,
from the relations below
$$\left\{\begin{array}{l}
(m^+)^{(i,2)}_{(i,2)}=K_i(m^+)^{(i+1,2)}_{(i+1,2)},i=1,2,\\
e^{(3,2)}(m^+)^{(i,2)}_{(i,2)}
=R_{\mathbb{C}^3}{}^{(i,3)}_{(i,3)}R_{(\mathbb{C}^2)^*}{}^{(2,2)}_{(2,2)}(m^+)^{(i,2)}_{(i,2)}e^{(3,2)},
\\
e^{(3,2)}(m^+)^{(i+1,2)}_{(i+1,2)}=R_{\mathbb{C}^3}{}^{(i+1,3)}_{(i+1,3)}R_{(\mathbb{C}^2)^*}{}^{(2,2)}_{(2,2)}(m^+)^{(i+1,2)}_{(i+1,2)}e^{(3,2)},
\end{array}
\right.
$$
we obtain
$e^{(3,2)}K_i=\frac{R_{\mathbb{C}^3}{}^{(i,3)}_{(i,3)}}{R_{\mathbb{C}^3}{}^{(i+1,3)}_{(i+1,3)}}K_ie^{3,2}=
\left\{
\begin{array}{ll}
K_ie^{3,2}& \mbox{when}~~i=1;\\
q^{-1}K_ie^{3,2}&\mbox{when}~~i=2.
\end{array}
\right.
$

On the other hand, from the relations below
$$\left\{\begin{array}{l}
(m^+)^{(3,1)}_{(3,1)}=K_4(m^+)^{(3,2)}_{(3,2)},\\
e^{(3,2)}(m^+)^{(3,2)}_{(3,2)}=R_{\mathbb{C}^3}{}^{(3,3)}_{(3,3)}R_{(\mathbb{C}^2)^*}{}^{(2,2)}_{(2,2)}(m^+)^{(3,2)}_{(3,2)}e^{(3,2)}\\
e^{(3,2)}(m^+)^{(3,1)}_{(3,1)}=R_{\mathbb{C}^3}{}^{(3,3)}_{(3,3)}R_{(\mathbb{C}^2)^*}{}^{(1,2)}_{(1,2)}(m^+)^{(3,1)}_{(3,1)}e^{(3,2)}
\end{array}
\right.
$$
we have
$e^{(3,2)}K_4=\frac{R_{(\mathbb{C}^2)^*}{}^{(1,2)}_{(1,2)}}{R_{(\mathbb{C}^2)^*}{}^{(2,2)}_{(2,2)}}K_4e^{(3,2)}=q^{\frac{1}{2}}K_4e^{(3,2)}.
$

\item
The cross relations between the new generator $e^{(3,2)}$ and $F_i,i=1,2,3,4$:

The relations $F_ie^{(3,2)}
=\frac{R_{\mathbb{C}^3}{}^{(3,i+1)}_{(3,i+1)}R_{(\mathbb{C}^2)^*}{}^{(2,2)}_{(2,2)}}
{R_{\mathbb{C}^3}{}^{(3,i+1)}_{(3,i+1)}R_{(\mathbb{C}^2)^*}{}^{(2,2)}_{(2,2)}}
e^{(3,2)}F_i=e^{(3,2)}F_i,i=1,2$ follow from 
$$
\left\{\begin{array}{l}
(m^-)^{(i+1,2)}_{(i,2)}=(q-q^-1)(m^-)^{(i+1,2)}_{(i+1,2)}F_i,\\
(m^-)^{(i+1,2)}_{(i,2)}e^{(3,2)}=R_{\mathbb{C}^3}{}^{(3,i+1)}_{(3,i+1)}R_{(\mathbb{C}^2)^*}{}^{(2,2)}_{(2,2)}e^{(3,2)}(m^-)^{(i+1,2)}_{(i,2)}
\\
(m^-)^{(i+1,2)}_{(i+1,2)}e^{(3,2)}=R_{\mathbb{C}^3}{}^{(3,i+1)}_{(3,i+1)}R_{(\mathbb{C}^2)^*}{}^{(2,2)}_{(2,2)}e^{(3,2)}(m^-)^{(i+1,2)}_{(i+1,2)}
\end{array}
\right.
$$
Similarly,
from the relations below
$$\left\{\begin{array}{l}
(m^-)^{(3,2)}_{(3,1)}=(q^{\frac{1}{2}}-q^{-\frac{1}{2}})(m^-)^{(3,2)}_{(3,2)}F_4,\\
(m^-)^{(3,2)}_{(3,1)}e^{(3,2)}=R_{\mathbb{C}^3}{}^{(3,3)}_{(3,3)}R_{(\mathbb{C}^2)^*}{}^{(2,2)}_{(2,2)}e^{(3,2)}(m^-)^{(3,2)}_{(3,1)}\\
(m^-)^{(3,2)}_{(3,2)}e^{(3,2)}=R_{\mathbb{C}^3}{}^{(3,3)}_{(3,3)}R_{(\mathbb{C}^2)^*}{}^{(2,2)}_{(2,2)}e^{(3,2)}(m^-)^{(3,2)}_{(3,2)},
\end{array}
\right.
$$
we have
$F_4e^{(3,2)}=\frac{R_{\mathbb{C}^3}{}^{(3,3)}_{(3,3)}R_{(\mathbb{C}^2)^*}{}^{(2,2)}_{(2,2)}}
{R_{\mathbb{C}^3}{}^{(3,3)}_{(3,3)}R_{(\mathbb{C}^2)^*}{}^{(2,2)}_{(2,2)}}
e^{(3,2)}F_4=e^{(3,2)}F_4$.
\item
The cross relations between the new ${K_3}$ and $E_1,E_2,E_4$:

According to the following cross relations in $U_q^{\text{ext}}(\mathfrak{sl}_3\oplus\mathfrak{sl}_2)\otimes
k[c,c^{-1}]$, namely,
\begin{gather*}
E_iK_i=q^2K_iE_i, \ (i=1,2), \ E_1K_2=q^{-1}K_2E_1, \  E_2K_1=q^{-1}K_1E_2, \   E_4K_4=qK_4E_4,
\end{gather*}
we figure out the relations below
\begin{gather*}
E_1\Big((m^+)^{(3,2)}_{(3,2)}c^{-1}\Big)=q^{-\frac{2}{3}}q^{\frac{2}{3}}\Big((m^+)^{(3,2)}_{(3,2)}c^{-1}\Big)E_1=\Big((m^+)^{(3,2)}_{(3,2)}c^{-1}\Big)E_1,
\\
E_2\Big((m^+)^{(3,2)}_{(3,2)}c^{-1}\Big)=q^{\frac{1}{3}}q^{-\frac{4}{3}}\Big((m^+)^{(3,2)}_{(3,2)}c^{-1}\Big)E_2=q^{-1}\Big((m^+)^{(3,2)}_{(3,2)}c^{-1}\Big)E_2,
\\
E_4\Big((m^+)^{(3,2)}_{(3,2)}c^{-1}\Big)=q^{-\frac{1}{2}}\Big((m^+)^{(3,2)}_{(3,2)}c^{-1}\Big)E_4.
\end{gather*}
\end{enumerate}

Second, we verify the $q$-Serre relations between $e^{(3,2)}$ and $E_i,i=1,2,4$.

With the following relations in hand
$$ \left\{
\begin{array}{l}(m^{+})^{(i,2)}_{(i+1,2)}=(q-q^{-1})E_i(m^{+})^{(i+1,2)}_{(i+1,2)},\quad i=1,2,\\
e^{(3,2)}(m^{+})^{(i,2)}_{(i+1,2)}=R_{\mathbb{C}^3}{}^{i,3}_{i,3}R_{(\mathbb{C}^2)^*}{}^{2,2}_{2,2}(m^{+})^{(i,2)}_{(i+1,2)}e^{(3,2)}+
R_{\mathbb{C}^3}{}^{i,3}_{i+1,2}R_{(\mathbb{C}^2)^*}{}^{2,2}_{2,2}(m^{+})^{(i+1,2)}_{(i+1,2)}e^{(2,2)},
\end{array}
\right.
$$
we obtain
\begin{gather}
e^{(3,2)}E_1=E_1e^{(3,2)},\quad\mbox{when~}i=1;\\
e^{(2,2)}=e^{(3,2)}E_2-q^{-1}E_2e^{(3,2)},\quad\mbox{when~}i=2.\label{serre1}
\end{gather}
And the relation between $E_2$ and $e^{(2,2)}$ is \begin{equation}\label{s11}e^{(2,2)}E_2=\frac{R_{\mathbb{C}^3}{}^{2,2}_{2,2}R_{(\mathbb{C}^2)^*}{}^{2,2}_{2,2}}
{R_{\mathbb{C}^3}{}^{3,2}_{3,2}R_{(\mathbb{C}^2)^*}{}^{2,2}_{2,2}}E_2e^{(2,2)}=qE_2e^{(2,2)},\end{equation}
which follows from the following equalities
$$ \left\{
\begin{array}{l}
e^{(2,2)}(m^{+})^{(2,2)}_{(3,2)}=R_{\mathbb{C}^3}{}^{2,2}_{2,2}R_{(\mathbb{C}^2)^*}{}^{2,2}_{2,2}(m^{+})^{(2,2)}_{(3,2)}e^{(2,2)},\\ e^{(2,2)}E_2(m^{+})^{(3,2)}_{(3,2)}=R_{\mathbb{C}^3}{}^{2,2}_{2,2}R_{(\mathbb{C}^2)^*}{}^{2,2}_{2,2}E_2(m^{+})^{(3,2)}_{(3,2)}e^{(2,2)},\\ e^{(2,2)}(m^{+})^{(3,2)}_{(3,2)}=R_{\mathbb{C}^3}{}^{3,2}_{3,2}R_{(\mathbb{C}^2)^*}{}^{2,2}_{2,2}E_2(m^{+})^{(3,2)}_{(3,2)}e^{(2,2)}.
\end{array}
\right.
$$

Substituting relation \eqref{serre1} into \eqref{s11},
we obtain
$$e^{(3,2)}E_2^2-(q+q^{-1})E_2e^{(3,2)}E_2+E_2^2e^{(3,2)}=0.$$
So we have the following $q$-Serre relation in the new quantum group $U$:
$$ E_3E_2^2-(q+q^{-1})E_2E_3E_2+E_2^2E_3=0.$$

Substituting relation \eqref{serre1} into the relation $(\blacktriangle)$ of $i=2,j=2$,
we have
$$[e^{(3,2)}]^3E_2-(q^{-1}{+}1{+}q)[e^{(3,2)}]^2E_2e^{(3,2)}+(q{+}q^{-1}{+}1)e^{(3,2)}E_2[e^{(3,2)}]^2-E_2[e^{(3,2)}]^3=0,$$
which is the $q$-Serre relation as required
$$
(E_{3})^{3}E_{2}-
\left[
\begin{array}{c}
3\\
1
\end{array}
\right]
_{q^{\frac{1}{2}}}
(E_{3})^{2}E_{2}E_{3}+
\left[
\begin{array}{c}
3\\
2
\end{array}
\right]
_{q^{\frac{1}{2}}}
E_{3}E_{2}(E_{3})^{2}-E_{2}(E_{3})^{3}=0.
$$

According to the first relation and the cross relations below
$$
\left\{
\begin{array}{l}
(m^{+})^{(3,1)}_{(3,2)}=(q^{\frac{1}{2}}-q^{-\frac{1}{2}})\,E_4\,(m^{+})^{(3,2)}_{(3,2)},\\
e^{(3,2)}\,(m^{+})^{(3,1)}_{(3,2)}=R_{\mathbb{C}^3}{}^{3,3}_{3,3}\,R_{(\mathbb{C}^2)^*}{}^{1,2}_{1,2}\,(m^{+})^{(3,1)}_{(3,2)}\,e^{(3,2)}+
R_{\mathbb{C}^3}{}^{3,3}_{3,3}\,R_{(\mathbb{C}^2)^*}{}^{1,2}_{2,1}\,(m^{+})^{(3,2)}_{(3,2)}\,e^{(3,1)},\\
e^{(3,2)}\,(m^{+})^{(3,2)}_{(3,2)}=R_{\mathbb{C}^3}{}^{3,3}_{3,3}\,R_{(\mathbb{C}^2)^*}{}^{2,2}_{2,2}\,(m^{+})^{(3,2)}_{(3,2)}\,e^{(3,2)},\\
e^{(3,1)}\,(m^{+})^{(3,2)}_{(3,2)}=R_{\mathbb{C}^3}{}^{3,3}_{3,3}\,R_{(\mathbb{C}^2)^*}{}^{2,1}_{2,1}\,(m^{+})^{(3,2)}_{(3,2)}\,e^{(3,1)},
\end{array}
\right.
$$
we obtain
\begin{equation}\label{serre2}
e^{(3,1)}=e^{(3,2)}E_4-q^{-\frac{1}{2}}E_4e^{(3,2)}.
\end{equation}
We thus obtain the $q$-Serre relation as required
$$
\left(e^{(3,2)}\right)^2E_4-(q^{\frac{1}{2}}+q^{-\frac{1}{2}})e^{(3,2)}E_4e^{(3,2)}+E_4\left(e^{(3,2)}\right)^2=0,
$$
by substituting $e^{(3,2)}e^{(3,1)}=q^{\frac{1}{2}}e^{(3,1)}e^{(3,2)}$ into equality \eqref{serre2}.

On the other hand, associated with
$$
e^{(3,1)}\,E_4=\frac{R_{\mathbb{C}^3}{}^{3,3}_{3,3}\,R_{(\mathbb{C}^2)^*}{}^{1,1}_{1,1}}
{R_{\mathbb{C}^3}{}^{3,3}_{3,3}\,R_{(\mathbb{C}^2)^*}{}^{2,1}_{2,1}}\,E_4\,e^{(3,1)}=q^{\frac{1}{2}}\,E_4\,e^{(3,1)},
$$
which follows from
$$
e^{(3,1)}(m^{+})^{(3,1)}_{(3,2)}=R_{\mathbb{C}^3}{}^{3,3}_{3,3}\,R_{(\mathbb{C}^2)^*}{}^{1,1}_{1,1}\,(m^{+})^{(3,1)}_{(3,2)}\,e^{(3,1)},
$$
we get another required $q$-Serre relation
$$
e^{(3,2)}(E_4)^2-(q^{\frac{1}{2}}+q^{-\frac{1}{2}})\,E_4\,e^{(3,2)}\,E_4+(E_4)^2\,e^{(3,2)}=0.
$$

With these in mind,
we prove that the identification is a Hopf algebra homomorphism.
It is straightforward to see that this is an isomorphism by virtue of the constructions in \cite{HH1} and \cite{HH2}.
Here, we can also give a vector space basis of $\tilde{V}(R',R)$ in a rough way,
similar to \textbf{Step 5} in Theorem \ref{typeA}.

According to the algebraic relations of $\tilde{V}(R',R)$, we obtain the following reduction systems in the spirit of \cite{Ber},
denoted by $S_1, S_2$ and $S_3$.
$$S_1=\left\{\sigma_1=(e^2e^1,q^{\frac{1}{2}}e^1e^2),
		\sigma_2=(e^4e^3,q^{\frac{1}{2}}e^3e^4),
		\sigma_3=(e^6e^5,q^{\frac{1}{2}}e^5e^6)
\right\},$$
	$$S_2=\left\{\begin{array}{l}
		\sigma_4=(e^4e^1,q^{\frac{1}{2}}e^3e^2-q^{-1}e^2e^3+q^{-\frac{1}{2}}e^1e^4)\\
		\sigma_5=(e^6e^3,q^{\frac{1}{2}}e^5e^4-q^{-1}e^4e^5+q^{-\frac{1}{2}}e^3e^6)\\
		\sigma_6=(e^6e^1,q^{\frac{1}{2}}e^5e^2-q^{-1}e^2e^5+q^{-\frac{1}{2}}e^1e^6)\end{array}
\right\},$$
$$	S_3=
\left\{\begin{array}{l}
	\sigma_7=(e^3e^3e^1,(1+q)e^3e^1e^3-qe^1e^3e^3),
		\sigma_8=(e^4e^4e^2,(1+q)e^4e^2e^4-qe^2e^4e^4),\\
		\sigma_9=(e^5e^5e^3,(1+q)e^5e^3e^5-qe^3e^5e^5),
		\sigma_{10}=(e^6e^6e^4,(1+q)e^6e^4e^6-qe^4e^6e^6),\\
		\sigma_{11}=(e^5e^5e^1,(1+q)e^5e^1e^5-qe^1e^5e^5),
		\sigma_{12}=(e^6e^6e^2,(1+q)e^6e^2e^6-qe^2e^6e^6)\end{array}
\right\}.
$$
It is convenient to treat $q^{\frac{1}{2}}e^3e^2-q^{-1}e^2e^3,\, q^{\frac{1}{2}}e^5e^4-q^{-1}e^4e^5, \,q^{\frac{1}{2}}e^5e^2-q^{-1}e^2e^5$, and $(1+q)e^{i+2}e^i-qe^ie^{i+2}$, for $i=1,2,3,4$, and $(1+q)e^{i+4}e^i-qe^ie^{i+4}$ for $i=1, 2$ as single combinations, respectively.

For convenience,
we denote these combinations by $\mathcal{M}_{ij}$ ($\mathcal{F}_{ij}$) in the following,
respectively.
Moreover,
we find that $e^{i+1}$ never occurs to the left of $e^i$, $i=1,3,5$,
the $q$-combination $(1+q)e^{i+2}e^i-qe^ie^{i+2},i=1,2,3,4$ never occurs to the right of $e^{i+2}$,
and $(1+q)e^{i+4}e^i-qe^ie^{i+4},i=1,2$ never occurs to the right of $e^{i+4}$ from the reduction system.
The overlap ambiguities of $S_i,i=1,2,3$ are
\begin{equation}\label{over1}(\sigma_2,\sigma_7,e^4,e^3,e^3e^1),(\sigma_3,\sigma_9,e^6,e^5,e^5e^3),(\sigma_3,\sigma_{11},e^6,e^5,e^5e^1);
\end{equation}
\begin{equation}\label{over3}(\sigma_8,\sigma_1,e^4e^4,e^2,e^1),(\sigma_{10},\sigma_2,e^6e^6,e^4,e^3),(\sigma_{12},\sigma_{1},e^6e^6,e^2,e^1);
\end{equation}
\begin{equation}\label{over4}(\sigma_5,\sigma_7,e^6,e^3,e^3e^1),(\sigma_{10},\sigma_{4},e^6e^6,e^4,e^1).\end{equation}
First, adjoining the following relations
$$e^3e^3\left(q^{\frac{1}{2}}e^3e^2{-}q^{-1}e^2e^3\right){=}(1+q)q^{{-}\frac{1}{2}}e^3\left(q^{\frac{1}{2}}e^3e^2{-}q^{-1}e^2e^3\right)e^3
{-}\left(q^{\frac{1}{2}}e^3e^2{-}q^{{-}1}e^2e^3\right)e^3e^3,$$
$$(1+q)q^{-\frac{1}{2}}e^5\left(q^{\frac{1}{2}}e^5e^4-q^{-1}e^4e^5\right)e^5
-\left(q^{\frac{1}{2}}e^5e^4-q^{-1}e^4e^5\right)e^5e^5=e^5e^5\left(q^{\frac{1}{2}}e^5e^4-q^{-1}e^4e^5\right),$$
$$(1+q)q^{-\frac{1}{2}}e^5\left(q^{\frac{1}{2}}e^5e^2-q^{-1}e^2e^5\right)e^5
-\left(q^{\frac{1}{2}}e^5e^2-q^{-1}e^2e^5\right)e^5e^5=e^5e^5\left(q^{\frac{1}{2}}e^5e^2-q^{-1}e^2e^5\right)$$
ensures that the ambiguities in \eqref{over1} are resolvable.

Moreover, the relations below
$$q^{-\frac{1}{2}}e^4\mathcal{M}_{23}e^2-\mathcal{F}_{24}\mathcal{M}_{23}=\mathcal{M}_{23}\mathcal{F}_{24}
-q^{-\frac{1}{2}}e^2\mathcal{M}_{23}e^4,$$
$$q^{-\frac{1}{2}}e^6\mathcal{M}_{25}e^2-\mathcal{F}_{26}\mathcal{M}_{23}=\mathcal{M}_{25}\mathcal{F}_{26}
-q^{-\frac{1}{2}}e^2\mathcal{M}_{25}e^6,$$
$$q^{-\frac{1}{2}}e^6\mathcal{M}_{45}e^4-\mathcal{F}_{46}\mathcal{M}_{45}=\mathcal{M}_{45}\mathcal{F}_{46}
-q^{-\frac{1}{2}}e^4\mathcal{M}_{45}e^6,$$
ensure that the ambiguities in \eqref{over3} are resolvable.

Consequently, $\mathcal{F}_{24}$ never occurs to the left of $\mathcal{M}_{23}$, $e^4$ never occurs to the left of $\mathcal{M}_{23}$,
$\mathcal{F}_{26}$ never occurs to the left of $\mathcal{M}_{25}$,
$e^6$ never occurs to the left of $\mathcal{M}_{25}$,
$\mathcal{F}_{46}$ never occurs to the left of $\mathcal{M}_{45}$,
and $e^6$ never occurs to the left of $\mathcal{M}_{45}$.

The following relations
$$e^5e^3\mathcal{M}_{23}=(1+q)q^{\frac{1}{2}}e^5\mathcal{M}_{23}e^3-e^3e^5\mathcal{M}_{23},$$
$$q^{\frac{1}{2}}(1+q)e^6e^6e^3e^2+q^{\frac{1}{2}}(1+q)e^6\mathcal{M}_{45}e^2
=-q\mathcal{M}_{45}e^2e^6-\mathcal{M}_{45}\mathcal{F}_{26}-q^{-\frac{1}{2}}(1+q)e^3\mathcal{F}_{26}e^6$$
ensure that the ambiguities in \eqref{over4} are resolvable. This means that $\mathcal{M}_{45}e^2$,$\mathcal{M}_{23}e^4$ and $e^5\mathcal{M}_{23}$ correspond to root vectors.
Moreover,
$\mathcal{M}_{45}e^2$ never occurs to the right of $e^5$,
$\mathcal{M}_{23}e^4$ never occurs to the right of $e^5$,
and $e^5\mathcal{M}_{23}$ never occurs to the left of $e^2$.
In addition,
$e^5\mathcal{M}_{45}e^2$,$e^5\mathcal{M}_{23}e^4$,$e^5\mathcal{M}_{23}e^2$ are also root vectors.

Therefore,
we obtain that the elements arranged in the order
\begin{gather*}
e^1, e^2, \mathcal{M}_{23}, \mathcal{F}_{13}, \mathcal{M}_{25}, e^3, e^5\mathcal{M}_{23}e^2, e^5\mathcal{M}_{23}, e^5\mathcal{M}_{23}e^4, \mathcal{F}_{24}, e^4,\\
\mathcal{F}_{15}, \mathcal{F}_{35}, \mathcal{M}_{45}, \mathcal{M}_{45}e^2, e^5\mathcal{M}_{45}e^2, e^5,
\mathcal{F}_{26},\mathcal{F}_{46}, e^6
\end{gather*}
form a convex PBW-basis generators of the braided group $V(R',R)$.

This completes the proof.
\end{proof}

\medskip\noindent
{\bf Acknowledgements.}
N. Hu is supported in part by the Fundamental and Interdisciplinary Disciplines Breakthrough Plan of the Ministry of Education of China (JYB2025XDXM112), the NSFC (Grant No. 12171155), and the Science and Technology Commission of Shanghai Municipality (Grant No. 22DZ2229014).

\medskip\noindent
{\bf Conflict of Interest:} \
The authors declare that they have no conflict of interest.

\end{document}